\documentclass[10pt]{amsart}
\usepackage{amsthm,amsfonts,amssymb,euscript}

\usepackage{xcolor}
\renewcommand{\epsilon}{\varepsilon}
\newcommand{\newsection}[1]
{\subsection{#1}\setcounter{theorem}{0} \setcounter{equation}{0}
\par\noindent}

\newtheorem{theorem}{Theorem}

  \newtheorem{proposition}[theorem]{Proposition}
  \newtheorem{lemma}[theorem]{Lemma}
  \newtheorem{corollary}[theorem]{Corollary}

\def\bphi{{\boldsymbol{\phi} }}
\def\bpsi{{\boldsymbol{\psi} }}
\def\F{{\mathcal F}}
\def\R{{\mathbb{R}}}

\def\Q{{\mathcal Q}}

\def\E{{\mathcal E}}

\def\ch{\mbox{ch} (k)}
\def\chbb{\mbox{chb} (k)}
\def\trigh{\mbox{trigh} (k)}

\def\sh{\mbox{sh} (k)}
\def\shh{\mbox{sh}}
\def\chh{\mbox{ch}}
\def\shhb{\overline {\mbox{sh}}}
\def\th{\mbox{th} (k)}
\def\chhb{\overline {\mbox{  ch}}}
\def\shbb{\mbox{shb} (k)}

\def\shb{\overline{\mbox{sh}(k)}}
\def\z{\zeta}
\def\zb{\overline{\zeta}}
\def\FF{{\mathbb{F}}}
\def\L{\Lambda}

\theoremstyle{remark}
  \newtheorem{remark}[theorem]{Remark}
  
\theoremstyle{definition}

\newcommand{\supp}{\operatorname{supp}}
\newcommand{\tr}{\operatorname{tr}}

\newcommand{\I}{\text{i}}
\newcommand{\e}{\varepsilon}
\newcommand{\dx}{\langle \nabla_x\rangle^{\frac12}}
\newcommand{\dy}{\langle \nabla_y\rangle^{\frac12}}
\newcommand{\dxh}{| \nabla_x|^{\frac12}}
\newcommand{\dyh}{| \nabla_y|^{\frac12}}
\newcommand{\dxx}{\langle \nabla_{x+y}\rangle^{\frac12}}
\newcommand{\dxxh}{| \nabla_{x+y}|^{\frac12}}
\newcommand{\dyy}{\langle \nabla_{x-y}\rangle^{\frac12}}

\newcommand{\dt}{|\partial_{t}|^{\frac14}}
\newcommand{\lc}{low\,\, collapsing}

\newcommand{\ls}{\lesssim}
\newcommand{\gp}{\Gamma_p}
\newcommand{\gpt}{\Gamma_{p,T}}
\newcommand{\gc}{\Gamma_c}
\newcommand{\gct}{\Gamma_{c,T}}
\newcommand{\lpp}{\Lambda_p}
\newcommand{\lppt}{\Lambda_{p,T}}
\newcommand{\lcc}{\Lambda_c}
\newcommand{\lcct}{\Lambda_{c,T}}
\newcommand{\sr}{\mathcal{S}'_r}
 \newcommand{\bgd}{\begin{equation}}
 \newcommand{\ad}{\begin{aligned}}
 \newcommand{\ead}{\end{aligned}}
 \newcommand{\ed}{\end{equation}}

\begin{document}
\title[Global estimates for HFB]
{GLOBAL UNIFORM IN N ESTIMATES FOR
SOLUTIONS OF A SYSTEM OF
HARTREE--FOCK--BOGOLIUBOV TYPE IN THE GROSS--PITAVESKII REGIME.}
%
%
%
%

\keywords{Hartree--Fock--Bogoliubov, Gross--Pitaevskii}
\subjclass[2010]{35Q55, 81}

\author[]{Xiaoqi Huang}
\address[X.H.]{Department of Mathematics, University of Maryland, College Park. MD 20742}
\email{xhuang49@umd.edu}

\begin{abstract}
We extend the recent work of Chong et al., 2022 \cite{1} to the critical case. More precisely, we prove global in time, uniform in $N$ estimates 
for the solutions $\phi$, $\Lambda$ and $\Gamma$  of a coupled system of Hartree--Fock--Bogoliubov type with interaction potential $\frac1NV_N(x-y)=N^{2}v(N(x-y))$. We assume that the potential $v$ is small which satisfies some technical conditions, and the initial conditions have finite energy. The main ingredient is a sharp estimate for the linear 
Schr\"odinger equation with potential in 6+1 dimension, which may be of interest in its own right.
\end{abstract}

\maketitle

\def\vac{{\big\vert 0\big>}}
\def\fpsi{{\big\vert\psi\big>}}
\def\bphi{{\boldsymbol{\phi} }}
\def\bpsi{{\boldsymbol{\psi} }}
\def\F{{\mathcal F}}
\def\R{{\mathbb{R}}}
\def\I{{\mathcal I}}
\def\Q{{\mathcal Q}}
\def\ch{\mbox{\rm ch} (k)}
\def\cht{\mbox{\rm ch} (2k)}
\def\chbb{\mbox{\rm chb} (k)}
\def\trigh{\mbox{\rm trigh} (k)}
\def\p{\mbox{p} (k)}
\def\sh{\mbox{\rm sh} (k)}
\def\sht{\mbox{\rm sh} (2k)}
\def\shh{\mbox{\rm sh}}
\def\chh{\mbox{ \rm ch}}
\def\th{\mbox{\rm th} (k)}
\def\thb{\overline{\mbox{\rm th} (k)}}
\def\shhb{\overline {\mbox{\rm sh}}}
\def\chhb{\overline {\mbox{\rm ch}}}
\def\shbb{\mbox{\rm shb} (k)}
\def\chbt{\overline{\mbox{\rm ch}  (2k)}}
\def\shb{\overline{\mbox{\rm sh}(k)}}
\def\shbt{\overline{\mbox{\rm sh}(2k)}}
\def\z{\zeta}
\def\zb{\overline{\zeta}}
\def\FF{{\mathbb{F}}}
\def\S{{\bf S}}
\def\Sr{{\bf S}_{red}}
\def\W{{\bf W}}
\def\WW{{\mathcal W}}
\def\L{{\mathcal L}}
\def\E{{\mathcal E}}
\def\X{{\tilde X}}
\def\Nor{{\rm Nor}}

\def \kb{{\bf k}}
\def\x{{\bf x}}
\def\fhat{{\widehat{f}}}
\def \ghat{{\widehat{g}}}
\def \fhatb{{\widehat{f}^{\ast}}}
\def \ghatb{{\widehat{g}^{\ast}}}
\def \Sb{\mathbb S}
\def\Zb {\mathbb Z}
\def\Lm{\Lambda}
\def \Lmb{\overline{\Lambda}}
\def \Gam{\Gamma}
\def\Gamb{\overline{\Gamma}}
\def\Lmh{\widehat{\Lambda}}
\def\Fhat{\widehat{F}}
\def\Fhatb{\overline{\widehat{F}}}
\def\E{{\mathcal E}}
\def\Sop{{\bf S}^{t}_{x_{1},x_{2}}}
\def\Wop{{\bf S}^{\pm,t}_{x_{1},x_{2}}}
\def\Sp{{\bf S}^{t}_{x_{1}}}
\def\nb{\nabla}
\def\Dl{\Delta}
\def\Acal{{\mathcal A}}
\def\ub{\overline{u}}
\def\tr{{\rm tr}}
\def\dig{{\rm diag}}

\def\Hbb{{\mathbb H}}
\def\Ubb{{\mathbb U}}
\def\Nbb{{\mathbb N}}

\def\zvec{\vec{z}}
\def\xvec{\vec{x}}
\def\zvecp{\vec{z}^{\prime}}

\newsection{Introduction}

Consider the $N$--body linear Schr\"odinger equation which governs the time-evolution of $N$ boson systems
\bgd\label{111}
\ad
\Big(\frac1i\frac{\partial}{\partial_t}-\sum_{j=1}^N \Delta_{x_j}+\frac1N  V_N(x_i-x_j)\Big)\psi_N(t,\cdot)=0,
\ead\ed
where $x_i\in\R^3$, $N$ is large and $V_N(x)=N^{3}v(N x)$  for some Schwarz class potential $v$. The conditions on the potential will be discussed below. 
A physically appealing case concerns initial data forming a tensor product of the same one--particle state, in spirit of the Bose-Einstein condensation. We refer to \cite{lieb2005mathematics}  for extensive background on Bose-Einstein condensation.

The goal is to find a rigorous, simple approximation to $\psi_N$ which is consistent with
\bgd
\psi_{approx}(t,x_1,\dots,x_N)\sim \phi(t,x_1)\phi(t,x_2)\dots\phi(t,x_N)
\ed
in an appropriate sense, where $\phi$ is often called the mean--field limit.

In the stationary case,  a survey of results concerning the ground
state properties of the dilute bosonic gases can be found in \cite{lieb2005mathematics}.
In the time dependent case, in the work of Erd\"os, Schlein and Yau  \cite{erdHos2007rigorous, erdHos2009rigorous, erdHos2010derivation}, by using the BBGKY hierarchies and the density matrix $\gamma_{N,t}$ formalism, the convergence 
of the exact dynamics to the mean--field limit is asserted in the trace norm as $N\rightarrow \infty$, provided that the 
the mean--field limit $\phi(t,x)$ satisfies the Gross--Pitaevskii equation
\bgd
\frac1i\frac{\partial}{\partial_t}\phi-\Delta\phi+8\pi a_0 |\phi|^2\phi=0.
\ed
Here $a_0$ is the scattering length of the potential $v$, physically, the scattering length measures the effective range of
the potential $V$, see e.g., \cite{erdHos2010derivation} for a precise definition of scattering length. Also in the recent work of Pickl \cite{pickl2015derivation}, the Gross--Pitaevskii equation is derived using a different method.
We refer to reader to \cite{benedikter2016effective,benedikter2015quantitative,rodnianski2009quantum,elgart2006gross, boccato2019bogoliubov} for  more backgrounds on the problem of approximating the many-body Schrödinger dynamics and related results in the Gross--Pitaveskii regime.

 The Hartree-Fock-Bogoliubov type equations are derived in Fock space,
which has been used in order to get estimates
for the rate of convergence of the approximation to the exact solution. In the mean--field regime, we also refer to the earlier work of Rodnianski and Schlein \cite{rodnianski2009quantum}, which gives the convergence rate  for marginal densities $\gamma_{N,t}^{(k)}$, and the work of  Lewin, Nam and Schlein \cite{lewin2015fluctuations}, which gives efficient direct estimates in $L^2(\R^N)$ using the Hartree state approach. In the Gross--Pitaveskii regime, we refer to the work of Brennecke and Schlein \cite{brennecke2019gross} for explicit rate of condensation.

The ideas behind the use of HFB system is to introduce pair excitations(given by a Bogoliubov transformation 
$e^{\mathcal{B}}$) as a second order correction to the mean field approximation.
The exact evolution of \eqref{111} is approximated by a construction involving two function: the condensate $\phi(t,x)$ and a pair 
excitation function $k(t,x,y)$, and it is 
\bgd \psi_{approx}:= e^{-\sqrt{N}\mathcal{A}(\phi(t))}e^{\mathcal{B}(k(t))}\Omega
\ed
where
\bgd
\mathcal{A}(\phi)=\int dx \{\bar\phi(x)a_x-\phi(x)a^*_x\}
\ed
with $a^*_x$ and $a_x$ being the creation and annihilation operator in Fock space, and $e^{-\sqrt{N}\mathcal{A}(\phi)}$ is a unitary operator on Fock space, the Weyl operator. And 
\bgd
\mathcal{B}(k)=\frac12\int dxdy \{\bar k(t,x,y)a_xa_y-k(t,x,y)a^*_xa^*_y \}.\ed
The unitary operator $e^{\mathcal{B}(k)}$ is the representation of an real symplectric matrix. Also, $\Omega$ is the vacuum. See for instance \cite{5} for background on this construction.

In order for $\psi_{approx}$ to be an approximation to the exact evolution, $\phi$ and $k$ must satisfy certain PDEs. They were introduced in \cite{5} and independently in a different 
context in \cite{bach2016time}, see also \cite{grillakis2010second}, \cite{3} for more details on the  argument motivating the use of Hartree-Fock-Bogoliubov type equations. We also mention that in Benedikter, de Oliveira
and Schlein \cite{benedikter2015quantitative}, a similar approach is considered, where the authors impose the Gross--Pitaevskii equation equation for $\phi$ and define $k$ by an
explicit formula, and give a rate of convergence result in terms of the marginal density.

In the remaining of  this paper, we shall focus on the analysis of the systems of PDEs. In particular, our main result, which we used in the analysis of nonlinear system, is the estimate for linear Schr\"odinger equation stated in the next section. The linear estimate in section 2 generalizes the recent work of  Chong et al., 2022 \cite{1} to the critical case.

The functions described by these equations are:
the condensate $\phi (t,x)$ and the density matrices
\begin{equation}\label{gammadefin}
\Gamma(t,x_1,x_2)=\frac1N\Big(\overline\sh\circ\sh\Big)(t,x_1,x_2)+\bar\phi(t,x_1)\phi(t,x_2)
\end{equation}
\begin{equation}\label{lammadefin}
\Lm(t,x_1,x_2)=\frac{1}{2N}\sht(t,x_1,x_2)+\phi(t,x_1)\phi(t,x_2)
\end{equation}
where 
\begin{equation}
\begin{aligned}
\sh&=k+\frac{1}{3!}k\circ \bar k\circ k+\dots \\
\ch&=\delta(x-y)+\frac{1}{2!}\bar k\circ k+\dots
\end{aligned}
\end{equation}
Here $(u\circ v) (x,y)=\int u(x,z)v(z,y)dz$, the pair excitation function $k$ is an auxiliary function, which does not explicitly appear in the system.

There are several equivalent ways of expressing the equations, in this section we shall use a compact, matrix formulation as in \cite{1}. We separate the condensate part from
the pair interaction part: define $\gc=\bar\phi\otimes \phi$, $\lcc=\phi\otimes \phi$, $\gp=\frac1N\overline\sh\circ\sh$ and $\lpp=\frac{1}{2N}\sht$. Also denote $\rho(t,x)=\Gamma(t,x,x)$

To write the Hartree-Fock-Bogoliubov equations in matrix notation, denote $$V_N(x-y)=N^3v(N(x-y))$$ for some Schwarz class potential which will be discussed below. Define
$$\Omega=\begin{pmatrix}
-\Gamma & -\bar\Lm\\
\Lm &  \Gamma
\end{pmatrix}=\Psi+\Phi
$$
where
$$\Psi=\begin{pmatrix}
-\Gamma_p & -\bar\Lm_p\\
\Lm_p &  \Gamma_p
\end{pmatrix}
$$
$$\Phi=\begin{pmatrix}
-\Gamma_c & -\bar\Lm_c\\
\Lm_c &  \Gamma_c
\end{pmatrix}
$$
Finally, let
$$S_3=\begin{pmatrix}
-I & 0\\
0 & I
\end{pmatrix}
$$
where $I$ is the identity operator.

The evolution equations for $\Psi$ and $\Phi$ are 
\begin{align}\label{Psi}
\frac1i &\partial_t \Phi-[\Delta_x\delta(x-y)S_3, \Phi] \nonumber\\
&=-[\big(V_N*\rho(t,x)\big)\delta(x-y)S_3, \Phi]-[V_N\Psi^*, \Phi]
\end{align}
\begin{align}\label{Phi}
\frac1i &\partial_t \Psi-[\Delta_x\delta(x-y)S_3, \Psi] \nonumber\\
&=-[\big(V_N*\rho(t,x)\big)\delta(x-y)S_3, \Psi]-\frac{1}{2N}[S_3, V_N\Psi] \\
&\qquad\qquad-[V_N\Omega^*, \Psi]-\frac{1}{2N}[S_3, V_N\Phi] \nonumber
\end{align}

In addition, the condensate $\phi$ satisfies
\bgd
\ad
\Big\{&\frac1i\partial_t-\Delta_{x_1}\Big\}\phi(x_1)\\
&=-\int dy \{V_N(x_1-y)\Gamma(y,y)\}\phi(x_1)  \\
&\,\,\,\,-\int dy \{V_N(x_1-y)\Gamma_p(y,x_1)\}\phi(y)  \\
&\,\,\,\,+\int dy \{V_N(x_1-y)\lpp(x_1,y)\}\bar\phi(y) 
\ead
\ed
Here $A^*(x,y)=\bar A(y,x), [A,B]=A\circ B-B\circ A$ and $V_N$ acts as pointwise multiplication by $V_N(x-y)$. See \eqref{lambdap}-\eqref{gammac} for a scalar form of the above 
equations.
 
The solutions $\phi,\,\Lm_p,\,\Lm_c\,\gp$ and $\gc$ all depend on $N$. This has been suppressed to simplify the notation. However, we will always keep track
of dependence on $N$ in our estimates. 

Next we review the conserved quantities of these equations, see \cite{5} for more details. The first conserved quantities is the total
number of particles(normalized by division by N):
\begin{equation}\label{particle}
\tr\{ \Gamma(t)\}=\|\phi(t,\cdot)\|_{L^2(dx)}^2+ \frac1N\|\sh(t,\cdot,\cdot)\|_{L^2(dxdy)}^2=1
\end{equation}
From here we see that
\begin{equation}\label{lm 120}
 \|\Lm(t,\cdot,\cdot)\|_{L^2(dxdy)} \le C
\end{equation}

The second conserved quantity is the energy per particle
\begin{equation}\label{ene}
\begin{aligned}
E(t)=&\tr\{ \nabla_{x_1}\cdot\nabla_{x_2}\Gamma(t)\} \\
&+\frac12 \int dx_1dx_2\Big\{V_N(x_1-x_2)|\Lm(t,x_1,x_2)|^2\Big\} \\
&+\frac12 \int dx_1dx_2\Big\{V_N(x_1-x_2)|\Gamma(t,x_1,x_2)|^2\Big\} \\
&+\frac12 \int dx_1dx_2\Big\{V_N(x_1-x_2)\Gamma(t,x_1,x_1)\Gamma(t,x_2,x_2)\Big\} \\
&-\int dx_1dx_2\Big\{V_N(x_1-x_2)|\phi(t,x_1)|^2|\phi(t,x_2)|^2\Big\} 
\end{aligned}
\end{equation}
The above holds for any Schwarz potential $v$. In addition, in order to use the estimates of \cite{2}, we assume
\begin{equation}\label{vconditio}
\begin{aligned}
&v \,\,\,\text{is spherically symmetric and} \\
&v\ge 0, \frac{\partial v}{\partial r}(r)\le 0.
\end{aligned}
\end{equation}

For the initial conditions, we assume that 
\begin{equation}\label{initial condition}
\begin{aligned}
\tr \{\Gamma(0)\}\le C_0 \\
E(0)\le C_0\end{aligned}
\end{equation}
The size of the initial condition $C_0$ is important, as we shall see in the end of section 5, 
our assumption on the size of the potential $v$ depends on the size of initial data.  Thus, we keep track of the constant $C_0$ from now on.

Note that the kinetic energy is 
\begin{equation}\label{eneke}
\begin{aligned}
\tr\{& \nabla_{x_1}\cdot\nabla_{x_2}\Gamma(t)\}= \int dx\Big\{|\nabla_x\phi(t,x)|^2\Big\}  \\
&+\frac{1}{2N}\int dx_1dx_2 \Big\{|\nabla_{x_1}\sh(t,x_1,x_2)|^2+|\nabla_{x_2}\sh(t,x_1,x_2)|^2\Big\} 
\end{aligned}
\end{equation}

If we assume $E\le C_0$, then we have an $H^1$ estimate for $\Lm$, uniformly in time (and $N$):
\begin{equation}\label{eneke1}
\begin{aligned}
&\int dx_1dx_2 \Big\{|\nabla_{x_1}\Lm(t,x_1,x_2)|^2+|\nabla_{x_2}\Lm(t,x_1,x_2)|^2\Big\}\le  CC_0 \\
&\frac{1}{N}\int dx_1dx_2 \Big\{|\nabla_{x_1,x_2}\sht(t,x_1,x_2)|^2\Big\} \le  CC_0
\end{aligned}
\end{equation}
Also, $\Gamma$ satisfies the $H^2$ type estimate
\begin{equation}\label{lm 12}
 \||\nabla_{x_1}||\nabla_{x_2}|\Gamma(t)\|_{L^2(dxdy)} \le  CC_0
\end{equation}
See \cite{3}, \cite{5}, as well as \cite{bach2016time} for these conserved quantities.

By Plancherel theorem, we see  that \eqref{eneke1} implies that for all time $t$
\begin{equation}\label{lm 1233}
 \|\dx\dy\Lm_p(t)\|_{L^2(dxdy)} \le CC_0.
\end{equation}
as well as 
\begin{equation}\label{lm 12333}
 \|\dx\dy\Lm_c(t)\|_{L^2(dxdy)} \le C C_0. 
\end{equation}
Here, and also later on, $\dx$ means $(1-\Delta_x)^{\frac14}$, which is a Fourier multiplier with symbol $(1+|\xi|^2)^{\frac14}$, and similarly for $\dy$.

Similarly, by \eqref{eneke} and Plancherel, we also have 
\begin{equation}\label{lm 123'}
 \|\dx\dy\gp(t)\|_{L^2(dxdy)} \le C C_0.
\end{equation}
as well as 
\begin{equation}\label{lm 123''}
 \|\dx\dy\gc(t)\|_{L^2(dxdy)} \le  CC_0.
\end{equation}

In order to state the main result for this paper in the simplest possible form, we define the following partial Strichartz norms:
\bgd\label{restr0}
\ad
&\|\Lambda\|_{\mathcal{S}_{x, y}}\\
&=  \sup_{  (p, q) \, \, admissible}\|\Lambda\|_{L^{p}(dt) L^{q}(dx) L^2(dy)}\\
&\quad +\sup_{ ( p, q) \, \, admissible}\|\Lambda\|_{L^{p}(dt) L^{q}(dy) L^2(dx)}
\ead
\ed
Recall $(p,q)$ are admissible in $3+1$ dimensions if $\frac2p+\frac3q=\frac32$, $2\le p\le \infty$. 

The main result of this paper is 
\begin{theorem}\label{main1.0}
Let $\Lm=\Lm_p+\Lm_c$, $\Gamma=\gp+\gc$ be solutions of \eqref{Psi}, \eqref{Phi}, while the potential satisfies \eqref{vconditio} and \eqref{vcondition} below, and
the initial conditions satisfy \eqref{initial condition}. Then we have 
\begin{equation}\label{main eq}
\|\dx\dy\Lm\|_{\mathcal{S}_{x, y}}+\|\dx\dy\Gamma\|_{\mathcal{S}_{x, y}} \le C
\end{equation}
for some constant $C$ independent of $N$. The above estimate still hold if we replace $\Lm$ by $\lcc$ or $\lpp$, or replace $\Gamma$ by $\gc$ or $\gp$.
\end{theorem}
We also have a theorem for $\sht$ (without dividing by $N$).
\begin{theorem}\label{main2.0}
Let $\Lm$, $\Gamma$, $\phi$  be solutions of \eqref{Psi}, \eqref{Phi}, while the potential satisfies \eqref{vconditio} and \eqref{vcondition} below, and
the initial conditions satisfy \eqref{initial condition}. Assume also that 
$$ \|\sht(0,\cdot,\cdot)\|_{L^2}+ \|\overline\sh\circ \sh(0,\cdot,\cdot)\|_{L^2}\le C
$$
Then we have 
\begin{equation}\label{main eq1}
\|\sht\|_{\mathcal{S}_{x, y}}+\|\overline\sh\circ \sh\|_{\mathcal{S}_{x, y}} \le C
\end{equation}
\end{theorem}
\begin{remark}
Here \eqref{main eq1}  improves the results in \cite{1}, Theorem 1.3 in two ways. First, the potential $N^2v(N(x-y))$ we considered here represents a stronger interaction between particles, compared with the $N^{3\beta-1}v(N^\beta(x-y))$, $\beta<1$ type potentials considered in \cite{1}.
Second, the  $\|\sht\|_{\mathcal{S}_{x, y}} $ norm stays bounded uniformly in $N$,
compared with the $\log N$ growth in \cite{1}. Although our argument is written for the case where the potential is $N^2v(N(x-y))$, it also works for the $N^{3\beta-1}v(N^\beta(x-y))$ case. For example, the uniform in $N$ estimates in \eqref{main eq1} still hold for the case where the potential is equal to $N^{3\beta-1}v(N^\beta(x-y))$.

The above estimates also imply some estimates for $\sh$. In particular,
$$ \|\sh\|_{L^p(dx)L^2(dy)} \le C  \|\sht\|_{L^p(dx)L^2(dy)}
$$
This is because $\sh=\frac12 \sht \circ\ch^{-1}$ and $\ch^{-1}$ has bounded operator norm.
\end{remark}

Finally, we also have estimates for $\phi$. Define the standard Strichartz spaces
\begin{align*}
&\|\phi\|_{\mathcal{S}}\label{restr1}
=  \sup_{  (p,q) \, \, admissible}\|\phi\|_{L^{p}(dt) L^{q}(dx)}.
\end{align*}
\begin{corollary}\label{main3.0}
Under the assumptions of Theorem~\ref{main1.0}, and the additional assumption $\|\langle \nabla\rangle^{\frac12}\phi (0,\cdot)\|_{L^2}\le C$
, we have 
\begin{equation}\label{main eq2}
\|\langle \nabla\rangle^{\frac12}\phi \|_{\mathcal{S}}\le C.\end{equation}
\end{corollary}

We shall mention the difficulties surrounding equation \eqref{Psi} and \eqref{Phi}. Denote
$$\S=\frac1i\frac{\partial}{\partial t}-\Delta_x-\Delta_y
$$
$$\S_{\pm}=\frac1i\frac{\partial}{\partial t}-\Delta_x+\Delta_y
$$
Schematically, treating $V_N$ as $\delta$ and ignoring the constants, the equations become
\bgd\ad
\S&\lcc=\Gamma(t,x,x)\lcc(t,x,y)+\lpp(t,x,x)\gc(t,x,y)\\
\S&_\pm\gc=\Gamma(t,x,x)\gc(t,x,y)+\bar\Lm_p(t,x,x)\lcc(t,x,y)\\
\S&\lpp+\frac{V_N}{N}\lpp=\Gamma(t,x,x)\lpp(t,x,y)+\lpp(t,x,x)\gp(t,x,y)-\frac{V_N}{N}\lcc \\
&+\lcc(t,x,x)\gp(t,x,y) \\
\S&_\pm\gp=\Gamma(t,x,x)\gp(t,x,y)+\lpp(t,x,x)\lpp(t,x,y)+\bar\Lm_c(t,x,x)\lpp(t,x,y)
\ead\ed
Recall that $V_N(x)=N^{3}v(N x)$, and $\frac{V_N}{N}$ in the $\lpp$ equation satisfies the critical scaling in the sense that $\frac{V_N}{N}\in L^{3/2}$ uniformly in $N\ge 1$.
In the case of one body problem, there is a vast amount of literature in the study of Schr\"odinger operators $-\Delta+V$ with critically singular potentials from different aspects see e.g., \cite{yajima,rodnianski2001time,journe1991decay,bouclet2018uniform,huang2022uniform}. And $V\in L^{n/2}$ for dimension $n\ge 3$ is almost the minimal condition to ensure the Schr\"odinger operators
$-\Delta+V$ is bounded from below and self-adjoint.

For the two-body problem in our case, borrowing ideas in \cite{1}, our method for treating the nonlinear term requires the Strichartz-type estimates for $\langle \nabla_x\rangle^{\alpha}\langle \nabla_y\rangle^{\alpha}\Lm_{p \, \,\text{or}\,\, c}$ and $\langle \nabla_x\rangle^{\alpha}\langle \nabla_y\rangle^{\alpha}\Gamma_{p \, \,\text{or}\,\, c}$ with $\alpha\ge\frac12$ . However, the $\dx\dy$ derivative is the threshold  for the linear equation in the sense that if we apply $\langle \nabla_x\rangle^{\alpha}\langle \nabla_y\rangle^{\alpha}$ with $\alpha>\frac12$ to the terms $\frac{V_N}{N}\lpp$ or $\frac{V_N}{N}\lcc$ in the $\lpp$ equation, we get a singularity which is essentially $N^{2\alpha-1}V_N\Lm$.  And since $\|N^{2\alpha-1}V_N \|_{L^p}\rightarrow \infty$ as $N\rightarrow \infty$  for any $1\le p \le \infty$, such terms can not be treated using standard techniques.

The structure of the rest of the paper is the following. In section 2, we list the notations used in this paper and state our main estimate, Theorem~\ref{main} for the linear Schr\"odinger equation in $6+1$ dimension. In sections 3 and 4, we prove
Theorem~\ref{main}. In section 5, we prove Theorem~\ref{main1.0} using the linear estimate Theorem~\ref{main}. In section 6 and section 7, we use prove Theorem~\ref{main2.0} and Corollary~\ref{main3.0}. In the last section, we use Theorem~\ref{main} to prove a ``collapsing estimate'' for the linear equation, which may be of interest on its own right.
\bigskip

\noindent\textbf{Acknowledgement}. The author is indebted to Manoussos Grillakis and Matei Machedon for both suggesting the problem and for many discussions and comments which improved the exposition. The author would also like to thank Jacky. Chong and Zehua Zhao for comments on an earlier version of this paper.

\newsection{List of notations and statement of the main linear estimates.}

Let us define the partial Strichartz norms
\begin{align}
&\|\Lambda\|_{\mathcal{S}_{x, y}}\label{restr}\\
&=  \sup_{  (p,q) \, \, admissible}\|\Lambda\|_{L^{p}(dt) L^{q}(dx) L^2(dy)}\notag\\
&\quad +\sup_{  (p,q) \, \notag \, admissible}\|\Lambda\|_{L^{p}(dt) L^{q}(dy) L^2(dx)}.
\end{align}
where the pair $(p,q)$ is admissible in $3+1$ dimension if $\frac2p+\frac3q=\frac32$, $2\le p\le \infty$.

Define the full Strichartz norm
\begin{align}
&\|\Lambda\|_{\mathcal{S}}\label{restrx-y}\\
&=  \sup_{  (p,q) \, \, admissible
}\|\Lambda\|_{L^{p}(dt) L^{q}(dx) L^2(dy)}\notag\\
&\quad +\sup_{  (p,q) \, \notag \, admissible}\|\Lambda\|_{L^{p}(dt) L^{q}(dy) L^2(dx)}\\
&\quad +\sup_{  (p,q) \, \notag \, admissible}\|\Lambda\|_{L^{p}(dt) L^{q}(d(x-y)) L^2(d(x+y))}.
\end{align}
And define the restricted dual Strichartz norm, excluding the end-points $p'=2$, $p'=1$: let $p_1$ large and
$p_0> 2$ but close to $2$, for admissible pairs $(p, q)$, define
\begin{align}\label{add1}
\|G\|_{\mathcal S'_r}= \inf_{  p_1 \ge p \ge p_0} \{\|G\|_{L^{p'}(dt) L^{q'}(dx) L^2(dy)},
\|G\|_{L^{p'}(dt) L^{q'}(dy) L^2(dx)}\}.
\end{align}

Let us also recall the standard Littlewood-Paley decomposition. Let $\phi(x)$ such that $\hat\phi \in C_0^\infty$ and 
$\hat\phi(\xi)=1$ in $|\xi|<1$, $\hat\phi(\xi)=0$ in $|\xi|>2$. Define the $\phi_k$ for $k\ge 0$ by $\hat\phi_k(\xi)=\hat\phi(\frac{\xi}{2^k})$ and denote  
\begin{equation}\label{add2}
 P_{|\xi|<2^k}f=f*\phi_k
\end{equation}
so that the inverse Fourier transform of $\hat\phi(\frac{\xi}{2^k}) \hat{f}$ is $P_{|\xi|<2^k}f$.

Next let $\psi_0=\phi$ and define $\psi_k$ for $k\ge 1$ by $\hat\psi_k(\xi)=\hat\phi(\frac{\xi}{2^k})-\hat\phi(\frac{\xi}{2^{k-1}})$. We also denote
\begin{equation}\label{add3}
 P_{|\xi|\sim 2^k}f=f*\psi_k
\end{equation}

For later use, we shall also abuse our notation a bit and define, for an arbitrary positive constant $M$,
\bgd \label{2.8} P_{|\xi|< M}f=f*\phi(Mx),
\ed
for any fixed constant $M$. And define,
\bgd\label{2.9} P_{|\xi|\ge M}f=f- P_{|\xi|< M}f.
\ed
So for any two fixed constants $0<M<N$, 
$$ P_{N\le|\xi|< M}f=P_{|\xi|< M}f- P_{|\xi|< N}f.
$$

Now we can define the following two ``collapsing norms''. Let 
\begin{align}\label{full collapsing}
\|\Lambda\|_{collapsing} =  \big\|\Lambda\big\|_{L^{\infty}(d(x-y))L^2(dt) L^2(d(x+y))}.
\end{align}
And define
\begin{equation}\label{low collapsing}
\begin{aligned}
\|\Lambda\|_{low\,\, collapsing}=& \big\|P_{|\xi-\eta|<20N}\Lambda\big\|_{collapsing} +\big\|P_{|\xi|<20N}\Lambda\big\|_{collapsing} \\
&\qquad+ \big\|P_{|\eta|<20N}\Lambda\big\|_{collapsing}.
\end{aligned}
\end{equation}

If $A\ls B$, there is a constant C such that $A\le CB$, and we use $A\sim B$ to denote the case when $A\ls B$ and $B\ls A$.

Define 
$\dx f=(1-\Delta_x)^{\frac14}f$
such that the Fourier transform of $\dx f$ is $(1+|\xi|^2)^{\frac14} \hat f$, and similarly the Fourier transform of $\dy f$ is $(1+|\eta|^2)^{\frac14} \hat f$ for any $f\in L^2(\R^6)$.

Let $x,y\in \R^3$, define 
\bgd\ad
\S&=\frac1i\frac{\partial}{dt}-\Delta_x-\Delta_y \\
\S_{\pm}&=\frac1i\frac{\partial}{\partial t}-\Delta_x+\Delta_y
\ead\ed

Consider the equation
\begin{equation}\label{2.1}
\begin{aligned}
\S\Lambda(t,x,y)&=N^2v\big(N(x-y)\big)\Lambda(t,x,y)+G(t,x,y)\\
&\qquad+N^2v\big(N(x-y))H(t,x,y) \\
 \Lambda(0,\cdot)&=\Lambda_0
\end{aligned} 
\end{equation}
where we shall assume that  
\begin{equation}\label{vcondition}
v \,\,\text{is Schwarz}\,\,\,\, \sup_{1\le p\le \infty} \|v\|_p<\varepsilon, \,\, and \,\,\, \supp \hat{v}\subset B_1(0).
\end{equation}
Here $ \hat{v}$ denotes the Fourier transform and $B_1(0)$ denotes the unit ball in $\R^3$, and 
$\e$ is a fixed small constant to be specified later which depends on $C_0$ in \eqref{initial condition}. The condition $\supp \hat{v}\subset B_1(0)$  is essentially not required for 
our proof, but it can simplify our argument, see the discussion in the beginning of section 4 for more details.

The simplest form of theorem is 

\begin{theorem}\label{main}
Let $\Lambda$ satisfy \eqref{2.1}, and assume $v$ satisfy \eqref{vcondition}, we have
\begin{equation}\label{main26}
\begin{aligned}
\|&\dx\dy \Lm\|_{\mathcal{S}_{x, y}} + \|\dxx \Lm\|_{\lc} + \|\dt \Lm\|_{\lc} \\
& \ls \|\dx\dy G\|_{\mathcal{S}'_{r}} +\e \|\dx\dy H\|_{L^2(dt)L^6(x-y) L^2(d(x+y))} \\
&\quad+ \e  \|\dxx H\|_{collapsing} +\e  \|\dt H\|_{collapsing}\\
&\quad + \e  \|\dx H\|_{collapsing} +\e  \|\dy H\|_{collapsing}+\|\dx\dy \Lm_0\|_{L^2}.
\end{aligned}
\end{equation}
\end{theorem}
\begin{remark}
The main difficulty in proving the theorem is the presence of the term $N^2v\big(N(x-y)\big)\Lambda(t,x,y)$, where $N^2v\big(N(x-y)\in L^{3/2}$ satisfies the critical scaling.
The term $N^2v\big(N(x-y)\big)H(t,x,y)$ in \eqref{2.1} is a technical term which arises from the term $N^2v\big(N(x-y)\big)\Lambda(t,x,y)$, since in our application, we split 
$\Lm=\lpp+\lcc$ and take $\Lm=\lpp$ and $H=\lcc$ in \eqref{2.1}. And the presence of $H(t,x,y)$ does not lead to any essential difficulty in the proof of the above theorem.

We also remark that in the case $H=0$, by using the above theorem plus an abstract argument, one can replace the $\|\cdot\|_{\lc}$ norm on the left side 
of \eqref{main26}  by the $\|\cdot\|_{collapsing}$ norm and the same 
result still holds. See the discussion in section 8 for more details.
\end{remark}

All the implicit constants in $\ls$ are independent of $N$ and $\e$, and the choice of the small constant $\e$ in \eqref{vcondition} will depend on the implicit constants and $C_0$ in \eqref{initial condition}.

\newsection{Preliminary estimates for solutions to the linear Schr\"odinger equation}

We will use the following Strichartz estimate. In 6+1 dimensions,
\begin{theorem}[Theorem 2.4, 2.5 of \cite{2}]\label{2.1'}
Let $\S u=f+g, \,\,u(0,\cdot)=u_0$. Then
\begin{equation}\nonumber
\|u\|_{\mathcal{S}}\ls \|f\|_{L^2(dt)L^{\frac65}(x-y) L^2(d(x+y))}+\|g\|_{\sr}+\|u_0\|_{L^2}.
\end{equation}
\end{theorem}

Now we shall present the main theorem of this section.
 \begin{theorem} \label{Spsi}   Let
$\S u = f \,\,u(0,\cdot)=0$
Then
\begin{align*}
\|u\|_{\mathcal{S}_{x, y}} \lesssim
\big\|\langle\nabla_{x+y}\rangle^{\frac{1}{2}}f\big\|_{L^{1}(d(x-y))L^2(dt)L^2(d(x+y))}
+\big\|\dt f\big\|_{L^{1}(d(x-y))L^2(dt)L^2(d(x+y))}.
\end{align*}
\end{theorem}
\begin{remark}
A frequency localized version of the above Theorem appears in \cite{1}, Proposition 4.7, which is also a motivation of the above theorem. The ${\mathcal{S}_{x, y}}$ norm is crucial here, we do not expect the above estimate to be true if we replace it by the full Strichartz norm ${\mathcal{S}}$ as in \eqref{restrx-y}.
\end{remark}

\begin{proof}
The main ideas is to divide the frequency support of $u$ into several regions, and use Strichartz estimate for the regions where $\tau\sim |\xi|^2+|\eta|^2$, and use Sobolev for the remaining regions. As we shall see later,  the proof of Theorem~\ref{Spsi1}-\ref{Spsi2} below uses essentially the same idea.

To begin with,
we shall use the decomposition $u=\sum_{k=0}^\infty P_{|\xi-\eta|\sim 2^k}u$, where for the case $k=0$, we are abusing notations a bit by letting $P_{|\xi-\eta|\sim 1}u$ to denote the operator $P_{|\xi-\eta|<1}u$. We have  the square function estimate(see e.g., Lemma 3.5 in \cite{1}).
\begin{equation}\label{add4}
\begin{aligned}
\|u\|_{\mathcal{S}_{x, y}} &\sim  \|\Big(\sum_{k=0}^\infty |P_{|\xi-\eta|\sim 2^k}u|^2\Big)^{\frac12}\|_{\mathcal{S}_{x, y}} \\
&\lesssim\Big(\sum_{k=0}^\infty  \|P_{|\xi-\eta|\sim 2^k}u\|_{\mathcal{S}_{x, y}}^2\Big)^{\frac12}.
\end{aligned}
\end{equation}
We shall focus on the dyadic pieces where $k\ge 1$, since by the Strichartz estimate and the Sobolev estimate, one can easily show that
$$ \|P_{|\xi-\eta| <1}u\|_{\mathcal{S}_{x, y}} \ls \big\|\langle\nabla_{x+y}\rangle^{\frac{1}{2}}f\big\|_{L^{1}(d(x-y))L^2(dt)L^2(d(x+y))}.
$$

Now let $u_k=P_{|\xi-\eta|\sim 2^k}u$, $f_k=P_{|\xi-\eta|\sim 2^k}f$, and decompose $u_k=u_k^1+u_k^2+u_k^3$, where

\begin{equation}\label{123}
\begin{aligned}
 \S \, u_k^1&= P_{10|\tau|^{\frac12}\ge 2^k} f_k, \,\,\,\,\text{with initial conditions 0} \\
\mathcal{F} \, u_k^2&=\frac{\mathcal{F}\Big( P_{10|\tau|^{\frac12}\le 2^k} f_k\Big)}{\tau+|\xi|^2+|\eta|^2}, \,\,\text{this no longer has initial conditions 0} \\
 \S \, u_k^3&= 0, \,\,\,\,\text{a correction so that $u_k^2+u_k^3$ has initial condition 0}.
 \end{aligned}
 \end{equation}

For $u_k^1$, by the Strichartz estimate
\begin{align*}
\|u_k^1\|_{\mathcal{S}} &\lesssim \|f_k\|_{{L^2(dt)L^{\frac{6}{5}}(d(x-y))L^2(d(x+y))} }  \\
&\lesssim  \|f_k\|_{L^{\frac{6}{5}}(d(x-y))L^2(dt)L^2(d(x+y)) }  \\
&\lesssim 2^{\frac k2} \|f_k\|_{L^{1}(d(x-y))L^2(dt)L^2(d(x+y)) },
\end{align*}
where in the last line we used the fact that $f_k$ is frequency supported in $|\xi-\eta|\sim 2^k$ and Bernstein's inequality, which is a (elementary) generalization of the classical Bernstein's  inequality to $L^2$ valued functions.

Now we make another dyadic decomposition, write 
\begin{equation}\label{add5}
f_k=\sum_{\ell\ge 0} P_{|\tau|\sim2^{2k+\ell}} f_k=\sum_\ell f_{k, \ell}.
\end{equation}

Note that for each fixed $\ell, k$, we have 
\begin{align*}
\|f_{k,\ell}\|_{L^{1}(d(x-y))L^2(dt)L^2(d(x+y)) } \sim&  2^{-\frac{\ell}{4}-\frac k2}\|\dt f_{k,\ell}\|_{L^{1}(d(x-y))L^2(dt)L^2(d(x+y)) }  \\
\lesssim &  2^{-\frac{\ell}{4}-\frac k2}\|\dt P_{|\tau|\sim2^{2k+\ell}} f\|_{L^{1}(d(x-y))L^2(dt)L^2(d(x+y)) }  
\end{align*}
where in the first line we used Bernstein's  inequality, and in the second line we used the fact that
\begin{equation}\label{add6}
\|P_{|\xi-\eta|\sim2^{k}} f\|_{L^{1}(d(x-y))L^2(dt)L^2(d(x+y)) }\ls \| f\|_{L^{1}(d(x-y))L^2(dt)L^2(d(x+y))},
\end{equation}
which can be proved, for example, using the generalized Young's inequality on the space of $L^2$ valued functions.

Thus, by Minkowski's inequality
\begin{align*}
\Big(\sum_{k=0}^\infty  \|u_k^1\|_{\mathcal{S}_{x, y}}^2\Big)^{\frac12} &\lesssim \Big(\sum_{k=0}^\infty 2^{k}\|\sum_\ell f_{k, \ell}\|^2_{L^{1}(d(x-y))L^2(dt)L^2(d(x+y)) }  \Big)^{\frac12}   \\
&\lesssim\sum_\ell \Big(\sum_{k=0}^\infty 2^{k}\| f_{k, \ell}\|^2_{L^{1}(d(x-y))L^2(dt)L^2(d(x+y)) }  \Big)^{\frac12}   \\
&\lesssim \sum_{\ell}2^{-\frac{\ell}{4}}\Big(\sum_k \|\dt P_{|\tau|\sim2^{2k+\ell}} f\|^2_{L^{1}(d(x-y))L^2(dt)L^2(d(x+y)) }  \Big)^{\frac12} \\
&\lesssim \sum_{\ell}2^{-\frac{\ell}{4}}\Big(\|\big(\sum_k|\dt P_{|\tau|\sim2^{2k+\ell}} f|^2\big)^{\frac12}\|^2_{L^{1}(d(x-y))L^2(dt)L^2(d(x+y)) }  \Big)^{\frac12} \\
&\lesssim \|\dt f\|_{L^{1}(d(x-y))L^2(dt)L^2(d(x+y)) },
\end{align*}
where in the last line we used the square function estimate in $t$ variable.

\bigskip
For $u_k^2$, the denominator is comparable with $|\xi-\eta|^2+|\xi+\eta|^2 \ge 2^{2k}\ge 100\tau$. Thus, by Sobolev's estimates at an angle, which is Lemma 3.2 in \cite{1}, we have
\begin{equation}\label{add7}
\begin{aligned}
\|u_k^2\|&_{L^2(dt)L^6( dx) L^2(dy)}+\|u_k^2\|_{L^2(dt )L^6(dy) L^2(dx)} \\ 
&\lesssim\|\langle \nabla_{x+y}\rangle u_k^2\|_{L^2(dt)L^2( dx) L^2(dy)} \\
&\lesssim\|\Big(\langle \nabla_{x-y}\rangle+\langle \nabla_{x+y}\rangle\Big)^{-2}\langle \nabla_{x+y}\rangle f_k\|_{L^2(dt)L^2( dx) L^2(dy)}.
\end{aligned}
\end{equation}
This is also the place where we require the norm on the left side to be  ${\mathcal{S}_{x, y}}$, since we do not have Sobolev-type estimates like 
\begin{align*}
\|u_k^2\|_{L^2(dt )L^6(d(x-y)) L^2(d(x+y))} 
\lesssim\|\langle \nabla_{x+y}\rangle u_k^2\|_{L^2(dt)L^2( dx) L^2(dy)}.
\end{align*}

Let us first assume $|\xi+\eta|\le |\xi-\eta|$, and make the 
 decomposition 
 \begin{equation}\label{add8}
 f_k=\sum_{j=0}^k P_{|\xi+\eta|\sim 2^{k-j}}f_k=\sum_jf_{k,j}.
 \end{equation}

Then for each fixed $k, j$, we have

\begin{equation}\label{add9}
\begin{aligned}
\|\Big(&\langle \nabla_{x-y}\rangle+\langle \nabla_{x+y}\rangle\Big)^{-2}\langle \nabla_{x+y}\rangle f_{k,j}\|_{L^2(dt)L^2( dx) L^2(dy)} \\
 \sim&\, 2^{-\frac{j}{2}}\|\langle \nabla_{x-y}\rangle^{-\frac32} \langle \nabla_{x+y}\rangle^{\frac12} f_{k,j}\|_{L^2(dt)L^2( dx) L^2(dy)}  \\
\lesssim & \, 2^{-\frac{j}{2}}\| \langle \nabla_{x+y}\rangle^{\frac12} f_{k,j}\|_{L^{1}(d(x-y))L^2(dt)L^2(d(x+y)) } \\
\lesssim & \, 2^{-\frac{j}{2}}\| \langle \nabla_{x+y}\rangle^{\frac12} P_{|\xi+\eta|\sim 2^{k-j}}f\|_{L^{1}(d(x-y))L^2(dt)L^2(d(x+y)) } ,
\end{aligned}
\end{equation}
where in the third line we used Bernstein's  inequality, and in the last line we used the fact that
$$\|f_k\|_{L^{1}(d(x-y))L^2(dt)L^2(d(x+y)) }\ls \| f\|_{L^{1}(d(x-y))L^2(dt)L^2(d(x+y))}. 
$$

Thus, by Minkowski's inequality
\begin{align*}
\Big(\sum_{k=0}^\infty  \|u_k^2\|&^2_{L^2(dt)L^6( dx) L^2(dy)}+\|u_k^2\|_{L^2(dt )L^6(dy) L^2(dx)}^2\Big)^{\frac12} \\
 &\lesssim \Big(\sum_{k=0}^\infty \|\sum_j \big(\langle \nabla_{x-y}\rangle+\langle \nabla_{x+y}\rangle\big)^{-2}\langle \nabla_{x+y}\rangle f_{k,j}\|^2_{L^{2}(d(x-y))L^2(dt)L^2(d(x+y)) }  \Big)^{\frac12}   \\
&\lesssim\sum_j \Big(\sum_{k} \| \big(\langle \nabla_{x-y}\rangle+\langle \nabla_{x+y}\rangle\big)^{-2}\langle \nabla_{x+y}\rangle f_{k, \ell}\|^2_{L^{2}(d(x-y))L^2(dt)L^2(d(x+y)) }  \Big)^{\frac12}   \\
&\lesssim \sum_{j}2^{-\frac{j}{2}}\Big(\sum_k \| \langle \nabla_{x+y}\rangle^{\frac12} P_{|\xi+\eta|\sim 2^{k-j}}f\|^2_{L^{1}(d(x-y))L^2(dt)L^2(d(x+y)) }  \Big)^{\frac12} \\
&\lesssim \sum_{j}2^{-\frac{j}{2}}\Big(\|\big(\sum_k| \langle \nabla_{x+y}\rangle^{\frac12} P_{|\xi+\eta|\sim 2^{k-j}}f|^2\big)^{\frac12}\|^2_{L^{1}(d(x-y))L^2(dt)L^2(d(x+y)) }  \Big)^{\frac12} \\
&\lesssim \| \langle \nabla_{x+y}\rangle^{\frac12} f\|_{L^{1}(d(x-y))L^2(dt)L^2(d(x+y)) } 
\end{align*}
where in the last line we used the square function estimate in $x+y$ variable.

The case $|\xi+\eta|\ge |\xi-\eta|$ is similar.

For the other endpoint $p=\infty$, define
\begin{equation}\label{add10}
u_{k}^2=\sum_{0\le \ell\le k/2} P_{|\tau|\sim 2^{2k-\ell}} u_k^2=\sum_\ell u^2_{k, \ell}
\end{equation}
and similarly 
\begin{equation}\label{add11}
f_{k}=\sum_{0\le \ell\le k/2} P_{|\tau|\sim 2^{2k-\ell}} f_k=\sum_\ell f_{k, \ell}.
\end{equation}

Note that for each fixed $\ell, k$, we have 
\begin{equation}\label{add12}
\begin{aligned}
\| u^2_{k,\ell}\|_{L^\infty(dt)L^2( dx) L^2(dy)} \lesssim&  \|\int\Big|\frac{\mathcal{F}\big( f_{k,\ell}\big)}{\tau+|\xi|^2+|\eta|^2}\Big| d\tau\|_{L^2( d\xi) L^2(d\eta)}\\
\sim & 2^{k/2-\ell/4} \| \langle \nabla_{x-y}\rangle^{-2}\dt  f_{k,\ell}\|_{L^2(dt)L^2( dx) L^2(dy)} \\
 \sim&\,2^{-\ell/4}\|\langle \nabla_{x-y}\rangle^{-\frac32}\dt  f_{k,\ell}\|_{L^2(dt)L^2( dx) L^2(dy)}  \\
\lesssim & \,2^{-\ell/4}\| \dt f_{k,\ell}\|_{L^{1}(d(x-y))L^2(dt)L^2(d(x+y)) } \\
\lesssim & \,2^{-\ell/4}\| \dt P_{|\tau|\sim 2^{2k-\ell}} f \|_{L^{1}(d(x-y))L^2(dt)L^2(d(x+y)) },
\end{aligned}
\end{equation}
where in the third line we used Bernstein's  inequality, and in the last line we used the fact that
$$\|f_k\|_{L^{1}(d(x-y))L^2(dt)L^2(d(x+y)) }\ls \| f\|_{L^{1}(d(x-y))L^2(dt)L^2(d(x+y))}. 
$$

By Minkowski's inequality
\begin{align*}
\Big(\sum_{k=0}^\infty  \| u_{k}^2\|&^2_{L^\infty(dt)L^2( dx) L^2(dy)}\Big)^{\frac12} \\
&=\Big(\sum_{k=0}^\infty  \|\sum_\ell u_{k,\ell}^2\|^2_{L^\infty(dt)L^2( dx) L^2(dy)}\Big)^{\frac12}  \\
&\lesssim \sum_{\ell=0}^\infty \Big(\sum_{k: k\ge 2\ell} \| u_{k,\ell}^2\|^2_{L^\infty(dt)L^2( dx) L^2(dy)}\Big)^{\frac12}   \\
&\lesssim\sum_{\ell=0}^\infty 2^{-\ell/4}\Big(\sum_{k: k\ge 2\ell} \| \dt P_{|\tau|\sim 2^{2k-\ell}} f \|^2_{L^{1}(d(x-y))L^2(dt)L^2(d(x+y)) } \Big)^{\frac12}  \\
&\lesssim\|\dt f\|_{L^{1}(d(x-y))L^2(dt)L^2(d(x+y)) } 
\end{align*}
where in the last line we used the square function estimate in $t$ variable.

To deal with $u_k^3$, note that since $u_k^3$ is solution to free Schr\"odinger, 
\begin{equation}\label{add13}
\begin{aligned}
\|u_k^3\|_{\mathcal{S}_{x, y}} &\lesssim  \|u_k^3(0,\cdot,\cdot)\|_{L^2( dx) L^2(dy)} \\
&= \|u_k^2(0,\cdot,\cdot)\|_{L^2( dx) L^2(dy)} \\
&\le   \|u_k^2\|_{L^\infty(dt)L^2( dx) L^2(dy)}.
\end{aligned}
\end{equation}
Thus it can be treated as in the previous case.

 \end{proof}

 \begin{theorem} \label{Spsi1}   Let
$\S u = f \,\,u(0,\cdot)=0$.
Then
\begin{equation}\label{xy dt}
\begin{aligned}
\|\dyy u\|_{\mathcal{S}_{x, y}} \lesssim &
\big\|\langle\nabla_{x+y}\rangle^{\frac{1}{2}}f\big\|_{L^2(dt)L^{\frac65}(d(x-y))L^2(d(x+y))} \\
&+\big\|\dt f\big\|_{L^2(dt)L^{\frac65}(d(x-y))L^2(d(x+y))}.
\end{aligned}
\end{equation}
\end{theorem}

\begin{proof}

We shall use the decomposition $u=\sum_{k=0}^\infty P_{|\xi-\eta|\sim 2^k}u$, and  the square function estimate
\begin{align*}
\|u\|_{\mathcal{S}_{x, y}} &\sim  \|\Big(\sum_{k=0}^\infty |P_{|\xi-\eta|\sim 2^k}u|^2\Big)^{\frac12}\|_{\mathcal{S}_{x, y}} \\
&\lesssim\Big(\sum_{k=0}^\infty  \|P_{|\xi-\eta|\sim 2^k}u\|_{\mathcal{S}_{x, y}}^2\Big)^{\frac12}.
\end{align*}
Since the right side of \eqref{xy dt} only involves $L^\frac65$-norm in the $x-y$ direction, we can add up the dyadic pieces in the right side using the square function estimate, thus it suffices to prove the Theorem for a single dyadic piece where $|\xi-\eta|\sim 2^k$(which would simplify our argument compare with the previous theorem). Here again for the case $k=0$, we are abusing notations a bit by letting $P_{|\xi-\eta|\sim 1}u$ to denote the operator $P_{|\xi-\eta|<1}u$, and the case $k=0$ is easy to handle by just using Strichartz.

Now let $u_k=P_{|\xi-\eta|\sim 2^k}u$ and $f_k=P_{|\xi-\eta|\sim 2^k}f$. We shall use the same decomposition as in \eqref{123}, write $u_k=u_k^1+u_k^2+u_k^3$.

For $u_k^1$, by Strichartz and Bernstein's  inequality
\begin{equation}\label{add14}
\begin{aligned}
\|\dyy u_k^1\|_{\mathcal{S}} &\lesssim \|\dyy f_k\|_{{L^2(dt)L^{\frac{6}{5}}(d(x-y))L^2(d(x+y))} }  \\
&\lesssim    \|\dt f_k\|_{L^2(dt)L^{\frac{6}{5}}(d(x-y))L^2(d(x+y)) } .
\end{aligned}
\end{equation}

\bigskip
For $u_k^2$, the denominator is comparable with $|\xi-\eta|^2+|\xi+\eta|^2 \ge 2^{2k}\ge 100\tau$. Thus, by Sobolev estimates at an angle, we have
\begin{equation}\label{add15}
\begin{aligned}
\|\dyy u_k^2\|&_{L^2(dt)L^6( dx) L^2(dy)}+\|\dyy u_k^2\|_{L^2(dt )L^6(dy) L^2(dx)} \\ 
&\lesssim\|\dyy \langle \nabla_{x+y}\rangle u_k^2\|_{L^2(dt)L^2( dx) L^2(dy)} \\
&\lesssim\|\Big(\langle \nabla_{x-y}\rangle+\langle \nabla_{x+y}\rangle\Big)^{-1}\dxx f_k\|_{L^2(dt)L^2( dx) L^2(dy)}\\
&\lesssim\| \langle\nabla_{x+y}\rangle^{\frac{1}{2}}f_{k}\|_{L^2(dt)L^{\frac65}(d(x-y))L^2(d(x+y)) }.
\end{aligned}
\end{equation}
For the other endpoint $p=\infty$, since we are in the region $10|\tau|^{\frac12}\le 2^k$, we have
\begin{equation}\label{add16}
\begin{aligned}
\|\dyy u^2_{k}&\|_{L^\infty(dt)L^2( dx) L^2(dy)}\\ \lesssim& 2^{k/2} \|\int\Big|\frac{\mathcal{F}\Big( P_{10|\tau|^{\frac12}\le \cdot2^k} f_k\Big)}{\tau+|\xi|^2+|\eta|^2}\Big| d\tau\|_{L^2( d\xi) L^2(d\eta)}\\
\lesssim & \,2^{k}\| \dt \langle \nabla_{x-y}\rangle^{-2} f_{k}\|_{L^2(dt)L^{2}(d(x-y))L^2(d(x+y)) } \\
\lesssim & \,\| \dt f_{k}  \|_{L^2(dt)L^{\frac65}(d(x-y))L^2(d(x+y)) }.
\end{aligned}
\end{equation}

\bigskip

To deal with $u_k^3$, if we repeat the argument in \eqref{add13}, we have 
\begin{align*}
\|u_k^3\|_{\mathcal{S}_{x, y}} &\lesssim   \|u_k^2\|_{L^\infty(dt)L^2( dx) L^2(dy)} .
\end{align*}
Thus it can be treated as in the previous case.

 \end{proof}

\begin{theorem} \label{Spsi2}   Let
\begin{align*}
&\S u = f,\,\,u(0, \cdot)=u_0
\end{align*}
We have
\begin{equation}
\begin{aligned}\label{dtdx}
\|\dt &u\|_{L^2(dt)L^{6}(d(x-y))L^2(d(x+y))} \\ 
&\lesssim \min\big\{\big\|\dt f\big\|_{L^2(dt)L^{\frac65}(d(x-y))L^2(d(x+y))}, \\
&\qquad\qquad \big\|\big(\langle\nabla_{x-y}\rangle^{\frac{1}{2}}+\langle\nabla_{x+y}\rangle^{\frac{1}{2}}\big)f\big\|_{L^2(dt)L^{\frac65}(d(x-y))L^2(d(x+y))}, \\
&\,\,\qquad\qquad \big\|\big(\langle\nabla_{x-y}\rangle^{\frac{1}{2}}+\langle\nabla_{x+y}\rangle^{\frac{1}{2}}\big)f\big\|_{\sr} \big\} \\
&+ \big\|\big(\langle\nabla_{x-y}\rangle^{\frac{1}{2}}+\langle\nabla_{x+y}\rangle^{\frac{1}{2}}\big)u_0\big\|_{L^2}.
\end{aligned} 
\end{equation}
\end{theorem}

\begin{proof}
 
For simplicity, we shall only treat the case where $|\xi+\eta|<|\xi-\eta|$, the case $|\xi+\eta|\ge |\xi-\eta|$ is similar and to some extent simpler.
 We shall use the decomposition $u=\sum_{k=0}^\infty P_{|\xi-\eta|\sim 2^k}u$, by using the square function estimate,
It suffices to prove the Theorem for a single dyadic piece where $|\xi-\eta|\sim 2^k$. The argument below also works for $k<0$.

Now let $u_k=P_{|\xi-\eta|\sim 2^k}u$ and $f_k=P_{|\xi-\eta|\sim 2^k}f$. We decompose $u_k=u_k^1+u_k^2+u_k^3+u_k^4$, where

\begin{align*}
 \S \, u_k^1&= P_{|\tau|^{\frac12}\sim 2^k} f_k, \,\,\,\,\text{with initial conditions 0} \\
\mathcal{F} \, u_k^2&=\frac{\mathcal{F}\Big( P_{10|\tau|^{\frac12}\le 2^k} f_k+P_{|\tau|^{\frac12}\ge 10\cdot2^k} f_k\Big)}{\tau+|\xi|^2+|\eta|^2}, \,\,\,\,\text{this no longer has initial conditions 0} \\
 \S \, u_k^3&= 0, \,\,\,\,\text{a correction so that $u_k^2+u_k^3$ has initial condition 0} \\
 \S \, u_k^4&= 0,\,\, \,\,\, u_k^4(0,\cdot)=P_{|\xi-\eta|\sim 2^k}u_0.
\end{align*}
It is easy to handle $u_k^4$, since in the case where $|\xi+\eta|<|\xi-\eta|$,
\begin{equation}\label{add17}
\begin{aligned}
\|\dt u_k^4&\|_{L^2(dt)L^{6}(d(x-y))L^2(d(x+y))}\\
 &\sim  2^{k/2}  \|e^{it\Delta}P_{|\xi-\eta|\sim 2^k}u_0\|_{L^2(dt)L^{6}(d(x-y))L^2(d(x+y)) }. \\
 &\lesssim \big\|\langle\nabla_{x-y}\rangle^{\frac{1}{2}}P_{|\xi-\eta|\sim 2^k}u_0\big\|_{L^2}.
\end{aligned}
\end{equation}

For $u_k^1$, since for $f_k$, we have $|\xi-\eta|^2+|\xi+\eta|\sim 2^{2k}$ and $\tau^{\frac12}\sim 2^{k}$,
 it is straightforward to check that
 the dual variable $\tau$ to $t$ for $u_k^1$ is also supported where $|\tau|^{1/2}\sim 2^{k}$, 
 by Strichartz estimates 
 \begin{equation}\label{add18}
\begin{aligned}
\|\dt u_k^1\|&_{L^2(dt)L^{6}(d(x-y))L^2(d(x+y))} \\
&\lesssim 2^{k/2}\|u_k^1\|_{L^2(dt)L^{6}(d(x-y))L^2(d(x+y))} \\
&\lesssim 2^{k/2} \|P_{|\tau|^{\frac12}\sim 2^k}f_k\|_{L^2(dt)L^{\frac{6}{5}}(d(x-y))L^2(d(x+y)) }  \\
&\lesssim   \|\dt f_k\|_{L^2(dt)L^{\frac{6}{5}}(d(x-y))L^2(d(x+y)) }.
\end{aligned}
\end{equation}
The same argument also gives
\begin{align*}
\|\dt u_k^1\|_{L^2(dt)L^{6}(d(x-y))L^2(d(x+y))}\lesssim   \|\langle\nabla_{x-y}\rangle^{\frac{1}{2}}f_k\|_{L^2(dt)L^{\frac{6}{5}}(d(x-y))L^2(d(x+y)) },
\end{align*}
as well as 
\begin{align*}
\|\dt u_k^1\|_{L^2(dt)L^{6}(d(x-y))L^2(d(x+y))}\lesssim   \|\langle\nabla_{x-y}\rangle^{\frac{1}{2}}f_k\|_{\sr}.
\end{align*}

\bigskip
For $u_k^2$, if $10|\tau|^{\frac12}\le 2^k$, the denominator is comparable with $|\xi-\eta|^2+|\xi+\eta|^2 \sim 2^{2k}$. Thus, by Sobolev's inequality in the $x-y$ direction, we have
\begin{align*}
\|\dt u_k^2\|_{L^2(dt)L^{6}(d(x-y))L^2(d(x+y))} 
\lesssim   \| \dt f_k\|_{L^2(dt)L^{\frac{6}{5}}(d(x-y))L^2(d(x+y)) }.
\end{align*}
The same argument also gives
\begin{align*}
\|\dt u_k^2\|_{L^2(dt)L^{6}(d(x-y))L^2(d(x+y))}\lesssim   \|\langle\nabla_{x-y}\rangle^{\frac{1}{2}}f_k\|_{L^2(dt)L^{\frac{6}{5}}(d(x-y))L^2(d(x+y)) }.
\end{align*}
It remains to show that,  if $10|\tau|^{\frac12}\le 2^k$, 
\begin{align}\label{uk2}
\|\dt u_k^2\|_{L^2(dt)L^{6}(d(x-y))L^2(d(x+y))}\lesssim   \|\langle\nabla_{x-y}\rangle^{\frac{1}{2}}f_k\|_{\sr}.
\end{align}
By interpolation, it suffices to show that
\begin{align}\label{uk21}
\|\dt u_k^2\|_{L^2(dt)L^{6}(d(x-y))L^2(d(x+y))}\lesssim   \|\langle\nabla_{x-y}\rangle^{\frac{1}{2}}f_k\|_{L^2(dt)L^{\frac{6}{5}}(dx)L^2(dy) },
\end{align}
\begin{align}\label{uk22}
\|\dt u_k^2\|_{L^2(dt)L^{6}(d(x-y))L^2(d(x+y))}\lesssim   \|\langle\nabla_{x-y}\rangle^{\frac{1}{2}}f_k\|_{L^2(dt)L^{\frac{6}{5}}(dy)L^2(dx) },
\end{align}
and 
\begin{align}\label{uk23}
\|\dt u_k^2\|_{L^2(dt)L^{6}(d(x-y))L^2(d(x+y))}\lesssim   \|\langle\nabla_{x-y}\rangle^{\frac{1}{2}}f_k\|_{L^1(dt)L^{2}(dx)L^2(dy) },
\end{align}
which would be stronger than \eqref{uk2} since it includes two endpoint cases.

The estimates \eqref{uk21} and \eqref{uk22} follow directly by Sobolev's estimates at an angle.
To prove \eqref{uk23}, first by Sobolev in the $x-y$ direction, we have
\begin{equation}\label{add19}
\begin{aligned}
\|\dt u_k^2\|_{L^2(dt)L^{6}(d(x-y))L^2(d(x+y))}\lesssim &  2^{-k/2}\|f_k\|_{L^2(dt)L^{2}(d(x-y))L^2(d(x+y)) }.
\\ =&  2^{-k/2}\|f_k\|_{L^2(dt)L^{2}(dx)L^2(dy) }
\end{aligned}
\end{equation}
and now by Bernstein's  inequality in the $t$ direction and $x-y $ direction.
\bgd\ad 2^{-k/2}\|f_k\|_{L^2(dt)L^{2}(dx)L^2(dy) }& \lesssim  2^{k/2}\|f_k\|_{L^1(dt)L^{2}(dx)L^2(dy) }  \\
& \lesssim  \|\langle\nabla_{x-y}\rangle^{\frac{1}{2}}f_k\|_{L^1(dt)L^{2}(dx)L^2(dy) }.
\ead\ed

\bigskip

If $|\tau|^{\frac12}\ge 10\cdot2^k$, the denominator is comparable with $\tau$.
Thus, by Sobolev in the $x-y$ direction, we have
\begin{equation}\label{add20}
\begin{aligned}
\|\dt u_k^2\|&_{L^2(dt)L^{6}(d(x-y))L^2(d(x+y))} \\
&\lesssim   \||\partial_{t}|^{-\frac{3}{4}}\langle\nabla_{x-y}\rangle^{2}f_k\|_{L^2(dt)L^{\frac{6}{5}}(d(x-y))L^2(d(x+y)) }\\
& \lesssim   \|\dt f_k\|_{L^2(dt)L^{\frac{6}{5}}(d(x-y))L^2(d(x+y)) }.
\end{aligned}
\end{equation}
The same argument also gives,
\begin{align*}
\|\dt u_k^2\|_{L^2(dt)L^{6}(d(x-y))L^2(d(x+y))}\lesssim   \|\langle\nabla_{x-y}\rangle^{\frac{1}{2}}f_k\|_{L^2(dt)L^{\frac{6}{5}}(d(x-y))L^2(d(x+y)) },
\end{align*}
 when $|\tau|^{\frac12}\ge 10\cdot2^k$.
It remains to show that \eqref{uk2} holds if  $|\tau|^{\frac12}\ge 10\cdot2^k$, which would be a consequence of \eqref{uk21}-\eqref{uk23}. And as before, 
the estimates \eqref{uk21} and \eqref{uk22} in this case follow by Sobolev's estimates at an angle.

To prove \eqref{uk23} when $|\tau|^{\frac12}\ge 10\cdot2^k$, first by Sobolev in the $x-y$ direction, we have
\begin{align*}
\|\dt u_k^2\|_{L^2(dt)L^{6}(d(x-y))L^2(d(x+y))} 
&\lesssim   \|\langle\partial_{t}\rangle^{-\frac{3}{4}}\langle\nabla_{x-y}\rangle f_k\|_{L^2(dt)L^{2}(d(x-y))L^2(d(x+y)) }\\
&= \|\langle\partial_{t}\rangle^{-\frac{3}{4}}\langle\nabla_{x-y}\rangle f_k\|_{L^2(dt)L^{2}(dx)L^2(dy) }.
\end{align*}
Now by Minkowski's inequality for $dt$ integral,
\begin{equation}\label{add21}
\begin{aligned}
\|\langle\partial_{t}\rangle^{-\frac{3}{4}}\langle\nabla_{x-y}\rangle f_k\|&_{L^2(dt)L^{2}(dx)L^2(dy) } \\
&=
 \|\tau^{-\frac34}\Big|\int e^{-it\tau}\langle\nabla_{x-y}\rangle f_k(t, \cdot) dt\Big|\|_{L^2(d\tau)L^{2}(dx)L^2(dy) }  \\
&\lesssim 2^{-k/2}
 \|\langle\nabla_{x-y}\rangle f_k\|_{L^1(dt)L^{2}(dx)L^2(dy) }\\   &\lesssim  \|\langle\nabla_{x-y}\rangle^{\frac{1}{2}}f_k\|_{L^1(dt)L^{2}(dx)L^2(dy) }.
\end{aligned}\end{equation}

\bigskip

To deal with $u_k^3$, note that since $u_k^3$ is solution to free Schr\"odinger equation,  and since we are assuming $|\xi-\eta|^2+|\xi+\eta|^2 \sim 2^{2k}$
\begin{equation}\label{add22}
\begin{aligned}
\|\dt u_k^3\|_{L^2(dt)L^{6}(d(x-y))L^2(d(x+y))}  &\lesssim 2^{k/2} \|u_k^3(0,\cdot,\cdot)\|_{L^2( dx) L^2(dy)} \\
&=2^{k/2}  \|u_k^2(0,\cdot,\cdot)\|_{L^2( dx) L^2(dy)} \\
&\le 2^{k/2}   \|u_k^2\|_{L^\infty(dt)L^2( dx) L^2(dy)} .
\end{aligned}
\end{equation}
Thus, it suffices to control $ \|u_k^2\|_{L^\infty(dt)L^2( dx) L^2(dy)} $.

First, if $|\tau|^{\frac12}\ge 10\cdot2^k$
\begin{equation}\label{add23}
\begin{aligned}
\| u^2_{k}\|_{L^\infty(dt)L^2( dx) L^2(dy)}
\lesssim&  \|\int\Big|\frac{\mathcal{F}\Big( P_{|\tau|^{\frac12}\ge 10\cdot2^k} f_k\Big)}{\tau+|\xi|^2+|\eta|^2}\Big| d\tau\|_{L^2( d\xi) L^2(d\eta)}\\
 \lesssim & 2^{-\frac{3k}{2}}  \| \dt f_{k}\|_{L^2(dt)L^2( dx) L^2(dy)} \\
\lesssim & \,2^{-k/2}\| \dt \langle \nabla_{x-y}\rangle^{-1} f_{k}\|_{L^2(dt)L^{2}(d(x-y))L^2(d(x+y)) } \\
\lesssim & \,2^{-k/2}\| \dt f_{k}  \|_{L^2(dt)L^{\frac65}(d(x-y))L^2(d(x+y)) }.
\end{aligned}
\end{equation}
Similarly, 
\begin{equation}\label{add24}
\begin{aligned}
\| u^2_{k}\|_{L^\infty(dt)L^2( dx) L^2(dy)}
\lesssim&  \|\int\Big|\frac{\mathcal{F}\Big( P_{|\tau|^{\frac12}\ge 10\cdot2^k} f_k\Big)}{\tau+|\xi|^2+|\eta|^2}\Big| d\tau\|_{L^2( d\xi) L^2(d\eta)}\\
 \lesssim & 2^{-k}  \| f_{k}\|_{L^2(dt)L^2( dx) L^2(dy)} \\
\lesssim & \,2^{-k/2}\| \langle \nabla_{x-y}\rangle^{\frac12} f_{k}\|_{L^2(dt)L^{\frac65}(d(x-y))L^2(d(x+y)) }.
\end{aligned}
\end{equation}
It remains to show that 
\begin{align}\label{uk20}
\| u^2_{k}\|_{L^\infty(dt)L^2( dx) L^2(dy)}
\lesssim \,2^{-k/2}\| \langle \nabla_{x-y}\rangle^{\frac12} f_{k}\|_{\sr}.
\end{align}
In this case, we won't prove \eqref{uk20} by interpolation since we do not know if one can show that
\begin{align}\label{uk200}
\| u^2_{k}\|_{L^\infty(dt)L^2( dx) L^2(dy)}
\lesssim \,2^{-k/2}\| \langle \nabla_{x-y}\rangle^{\frac12} f_{k}\|_{L^1(dt)L^{2}(dx)L^2(dy) },
\end{align}
which is also the reason why we have the restricted norm ${\sr}$ in the statement of the Theorem.

Instead, we shall prove \eqref{uk20} by showing that
\begin{align}\label{uk2000}
\| u^2_{k}\|_{L^\infty(dt)L^2( dx) L^2(dy)}
\lesssim \,2^{-k/2}\| \langle \nabla_{x-y}\rangle^{\frac12} f_{k}\|_{L^{p'}(dt)L^{q'}(dx)L^2(dy) },
\end{align}
as well as 
\begin{align}\label{uk20000}
\| u^2_{k}\|_{L^\infty(dt)L^2( dx) L^2(dy)}
\lesssim \,2^{-k/2}\| \langle \nabla_{x-y}\rangle^{\frac12} f_{k}\|_{L^{p'}(dt)L^{q'}(dy)L^2(dx) },
\end{align}
for all admissible pairs $(p,q)$, with $\frac2p=\frac 32-\frac 3q,\,\, 2\le p<\infty$.

To prove \eqref{uk2000} when $|\tau|^{\frac12}\ge 10\cdot2^k$,  for  admissible pair $(p,q)$ with $p<\infty$,
\begin{equation}\label{add25}
\begin{aligned}
\| u^2_{k}\|_{L^\infty(dt)L^2( dx) L^2(dy)}
\lesssim&  \|\int\Big|\frac{\mathcal{F}\Big( P_{|\tau|^{\frac12}\ge 10\cdot2^k} f_k\Big)}{\tau+|\xi|^2+|\eta|^2}\Big| d\tau\|_{L^2( d\xi) L^2(d\eta)}\\
\lesssim& 2^{-\frac{2k}{p}} \|\mathcal{F}\Big( P_{|\tau|^{\frac12}\ge 10\cdot2^k }f_k\Big)\|_{L^2( d\xi) L^2(d\eta)L^p(d\tau)}\\
 \lesssim &  2^{-\frac{2k}{p}}  \| f_{k}\|_{L^{p'}(dt)L^2( dx) L^2(dy)} \\
  \lesssim &  2^{-\frac{2k}{p}} 2^{(\frac 32-\frac 3q)k} \| f_{k}\|_{L^{p'}(dt)L^{q'}( dx) L^2(dy)} \\
   \lesssim &   \| f_{k}\|_{L^{p'}(dt)L^{q'}( dx) L^2(dy)} \\
\lesssim & \,2^{-k/2}\| \langle \nabla_{x-y}\rangle^{\frac12} f_{k}\|_{L^{p'}(dt)L^{q'}( dx) L^2(dy)},
\end{aligned}
\end{equation}
where we used H\"older's inequality in the second line, the Hausdorff-Young inequality in the third line, and Bernstein's inequality at an angle in the fourth and last line.
The proof of \eqref{uk20000} is similar.

\bigskip

If $10|\tau|^{\frac12}\le 2^k$, 
\begin{equation}\label{add26}
\begin{aligned}
\| u^2_{k}\|_{L^\infty(dt)L^2( dx) L^2(dy)}
\lesssim&  \|\int\Big|\frac{\mathcal{F}\Big( P_{10|\tau|^{\frac12}\le \cdot2^k} f_k\Big)}{\tau+|\xi|^2+|\eta|^2}\Big| d\tau\|_{L^2( d\xi) L^2(d\eta)}\\
 \lesssim & 2^{-\frac{3k}{2}}  \| \dt f_{k}\|_{L^2(dt)L^2( dx) L^2(dy)} \\
\lesssim & \,2^{-k/2}\| \dt \langle \nabla_{x-y}\rangle^{-1} f_{k}\|_{L^2(dt)L^{2}(d(x-y))L^2(d(x+y)) } \\
\lesssim & \,2^{-k/2}\| \dt f_{k}  \|_{L^2(dt)L^{\frac65}(d(x-y))L^2(d(x+y)) }.
\end{aligned}
\end{equation}
The same argument also gives,
\begin{align*}
\| u^2_{k}\|_{L^\infty(dt)L^2( dx) L^2(dy)}
\lesssim  \,2^{-k/2}\| \langle \nabla_{x-y}\rangle^{\frac12}f_{k}  \|_{L^2(dt)L^{\frac65}(d(x-y))L^2(d(x+y)) }.
\end{align*}
It remains to show that 
\begin{align}\label{uk20'}
\| u^2_{k}\|_{L^\infty(dt)L^2( dx) L^2(dy)}
\lesssim \,2^{-k/2}\| \langle \nabla_{x-y}\rangle^{\frac12} f_{k}\|_{\sr},
\end{align}
for the case $10|\tau|^{\frac12}\le 2^k$, which would be a consequence of \eqref{uk21}-\eqref{uk23}.

To prove \eqref{uk21}, 
\begin{align*}
\| u^2_{k}\|_{L^\infty(dt)L^2( dx) L^2(dy)}
\lesssim&  \|\int\Big|\frac{\mathcal{F}\Big( P_{|\tau|^{\frac12}\ge 10\cdot2^k} f_k\Big)}{\tau+|\xi|^2+|\eta|^2}\Big| d\tau\|_{L^2( d\xi) L^2(d\eta)}\\
 \lesssim & 2^{-k}  \| f_{k}\|_{L^2(dt)L^2( dx) L^2(dy)} \\
\lesssim & \,2^{-k/2}\| \langle \nabla_{x-y}\rangle^{\frac12} f_{k}\|_{L^2(dt)L^{\frac65}(dx)L^2(dy) }.
\end{align*}
The proof of \eqref{uk22} is similar.

To prove \eqref{uk23}, we use Bernstein's inequality in the $t$ direction
\begin{equation}\label{add27}
\begin{aligned}
 2^{-k}\|f_k\|_{L^2(dt)L^{2}(dx)L^2(dy) } &\lesssim  \|f_k\|_{L^1(dt)L^{2}(dx)L^2(dy) }  \\
& \lesssim  2^{-k/2} \|\langle\nabla_{x-y}\rangle^{\frac{1}{2}}f_k\|_{L^1(dt)L^{2}(dx)L^2(dy) }.
\end{aligned}\end{equation}

\end{proof}

\bigskip 

Now we shall present several lemmas that involve the collapsing norm.

\begin{lemma}\label{xco}
If$\,\, \S u=g, u(0,\cdot)=u_0$. Then
\begin{align*}
\|u\|_{collapsing} \lesssim \min\{\|\dx g\|_{\sr}+ \|\dx u_0\|_{L^2}, \|\dy g\|_{\sr}+ \|\dy u_0\|_{L^2}\}.
\end{align*}
\end{lemma}

We record that the above implies
\begin{lemma}\label{xco1}
If$\,\, \S u=g, u(0,\cdot)=u_0$. Then
\begin{align*}
\|\dx u&\|_{collapsing}+\|\dy u\|_{collapsing} \\
& \lesssim \|\dx\dy g\|_{\sr}+ \|\dx\dy u_0\|_{L^2}.
\end{align*}
\end{lemma}

We will also need
\begin{lemma}\label{xco2}
If$\,\, \S u=g, u(0,\cdot)=u_0$. Then
\begin{align}\label{ss}
\|\dxx u\|_{collapsing} 
 \lesssim \|\dx\dy g\|_{\sr}+ \|\dx\dy u_0\|_{L^2}.
\end{align}
\end{lemma}

The proof of Lemma~\ref{xco} and Lemma~\ref{xco2} are similar, for simplicity, we shall only present the proof of Lemma~\ref{xco2} here. 
The proof essentially follows from ideas in Lemma 5.1, 5.3 in \cite{3}.

\begin{proof}
We shall first prove the homogeneous estimate, let $ \S u=0$, with $u(0,\cdot)=u_0$.
Our goal is to show 
\begin{align}\label{hoco}
\sup_{x-y}\|\dxxh u\|_{L^{2}(dt)L^2(d(x+y) } 
 \lesssim \|\dxh\dyh u_0\|_{L^2}.
\end{align}
This is stronger than desired, since $\|\dxh\dyh u_0\|_{L^2}\ls \|\dx\dy u_0\|_{L^2}$, and also by Lemma~\ref{xco}, we have
\begin{align*}
\sup_{x-y}\| u\|_{L^{2}(dt)L^2(d(x+y) } 
 \lesssim \|\dx u_0\|_{L^2}\ls   \|\dx\dy u_0\|_{L^2}.
\end{align*}

To prove \eqref{hoco}
,  let $\widetilde{\Lambda}$ denote the space-time Fourier. For fixed $x-y$,
doing Cauchy-Schwarz with measures,
\begin{equation}\label{add28}
\begin{aligned}
&|\dxxh\widetilde{\Lambda(t, x-y, x+y)}(\tau, \xi+\eta)|^2 \\
&\lesssim
 \int \delta(\tau-|\xi|^2 -|\eta|^2)\frac{|\xi+\eta|}{|\xi||\eta|}d (\xi-\eta) \\
 &\ls \int \delta(\tau-|\xi|^2 -|\eta|^2)|\widehat{\nabla_x^{\frac{1}{2}} \nabla_y^{\frac{1}{2}}\Lambda_0}(\xi, \eta)|^2 d (\xi-\eta).
\end{aligned}
\end{equation}
In order to prove the estimate, we must show
\begin{align*}
\sup_{\tau, \xi}  \int \delta(\tau-|\xi|^2 -|\eta|^2)\frac{|\xi+\eta|}{|\xi||\eta|}d (\xi-\eta)\lesssim 1.
\end{align*}
Without loss of generality, consider the region $|\xi| \le | \eta|$. If $|\xi| \sim | \eta|$,
$\frac{|\xi+\eta|}{|\xi||\eta|}\lesssim \frac{1}{|\xi-\eta|}  $ and the integral can be evaluated in polar coordinates.
If $|\xi| << | \eta|$ then $|\xi+\eta| \sim |\xi-\eta|$
Writing $\frac{|\xi+\eta|}{|\xi||\eta|} \lesssim \frac{1}{|\xi|} \lesssim \frac{1}{|\xi-\eta|\sqrt{1-\cos(\theta)}}$
  where $\theta$ is the angle between $\xi-\eta$ and $\xi+\eta$, we estimate
\begin{align}
\sup_{\tau} \int_0^{\pi}\int \delta(\tau - \rho^2)\frac{1}{\rho\sqrt{1-\cos(\theta)}} \rho^2d \rho \sin(\theta) d \theta \lesssim
1.
\end{align}
\end{proof}

The inhomogeneous estimate \eqref{ss} just follows from the homogeneous estimate \eqref{hoco} and the Christ--Kiselev lemma.
More precisely, let $T_1=e^{it(\Delta_x+\Delta_y)}$, so $T_1: L^2(\R^6) \rightarrow L^p(dt)L^q(dx)L^2(dy)$ and $$T_1^*:
L^{p'}(dt)W^{\alpha, q'}(dx)H^{\alpha}(dy) \rightarrow H^{\alpha}(dx) H^{\alpha}(dy). $$
Fix $x-y$ and let $T_2: H^{\alpha}(dx) H^{\alpha}(dy)\rightarrow L^2(dt)H^\alpha(d(x+y))$ be the operator
$f\rightarrow \big(e^{it(\Delta_x+\Delta_y)}f\big) (t, x-y, x+y)$. Then the inhomogeneous estimate follows by applying the Christ--Kiselev lemma to $T_2T_1^*$.

\begin{lemma}\label{xco3}
If$\,\, \S u=g, u(0,\cdot)=u_0$. Then
\begin{align}\label{sss}
\|\dt u\|_{collapsing} 
 \lesssim \|\dx\dy g\|_{\mathcal{S}^\prime_{r}}+ \|\dx\dy u_0\|_{L^2}.
\end{align}
\end{lemma}
\begin{proof}
For the homogeneous estimate, it follows from the same argument as above. However, we can not apply Christ--Kiselev lemma here to get inhomogeneous estimate since $\dt$ does not commute with $1_{[0,t]}$ when we write out the solution using Duhamel's formula. Let $\S u=g$, with $u(0,\cdot)=0$, it suffices to prove
\begin{align}\label{ssss}
\|\dt u\|_{collapsing} 
 \lesssim \|\dx\dy g\|_{\mathcal{S}^\prime_{r}}.
\end{align}
To prove this, we shall decompose the Fourier support $\tau$ and $|\xi+\eta|$ of $u$ into finitely many regions.

Case 1: $\tau^{\frac12}\le {10}(1+|\xi+\eta|)$.

In this case, we have 
$$ \|\dt u\|_{collapsing}  \ls \|\dxx u\|_{collapsing},$$
thus the desired estimates follows from Lemma~\ref{xco2}.

\bigskip

Case 2: If $|\tau|^{\frac{1}{2}} > 2 (|\xi|+ |\eta|)$.

Write $u=u^1+u^2$, where
\begin{equation}\label{123'}
\begin{aligned}
\mathcal{F} \, u^1&=\frac{\mathcal{F}f}{\tau+|\xi|^2+|\eta|^2}, \,\,\text{this no longer has initial conditions 0} \\
 \S \, u^2&= 0, \,\,\,\,\text{a correction so that $u^1+u^2$ has initial condition 0}.
 \end{aligned}
 \end{equation}
In this case, it suffices to control $u_1$ since $u_2$ is only supported where $|\tau|=|\xi|^2+|\eta|^2$. The goodness about $u_1$ is that it has the same Fourier support with $f$.
The strategy is based on
\begin{align*}
\|\dt u_1 \|&_{L^{\infty}(d(x-y)) L^2(d(x+y)dt)}
=\|\tau^{\frac14}\frac{\mathcal{F}f}{\tau+|\xi|^2+|\eta|^2} \|_{L^{\infty}(d(x-y))L^2(d \tau) d (\xi+\eta)}\\
&\lesssim \|\int|\tau|^{\frac14} \Big|\frac{\mathcal{F}\Big(f\Big)}{\tau+|\xi|^2+|\eta|^2}\Big|d(\xi-\eta)\|_{L^2(d \tau d(\xi+\eta))}.
\end{align*}
It suffices to show
\begin{equation}\label{uk1}
\begin{aligned}
\|\int|\tau|^{\frac14} \Big|\frac{\mathcal{F}f}{\tau+|\xi|^2+|\eta|^2}\Big|d(\xi-\eta)\|_{L^2(d \tau d(\xi+\eta))}
\\ \lesssim \,\| |\nabla_{x}|^{\frac12}|\nabla_{y}|^{\frac12} f\|_{L^{p'}(dt)L^{q'}(dx)L^2(dy) },
\end{aligned}
\end{equation}
as well as 
\begin{equation}\label{uk12}
\begin{aligned}
\|\int|\tau|^{\frac14} \Big|\frac{\mathcal{F}f}{\tau+|\xi|^2+|\eta|^2}\Big|d(\xi-\eta)\|_{L^2(d \tau d(\xi+\eta))}
\\ \lesssim \,\| |\nabla_{x}|^{\frac12}|\nabla_{y}|^{\frac12} f\|_{L^{p'}(dt)L^{q'}(dy)L^2(dx) },
\end{aligned}
\end{equation}
for all admissible pairs $(p,q)$, where $\frac2p=\frac 32-\frac 3q,\,\, 2\le p<\infty$.

For $2\le p< \infty$, by Cauchy-Schwarz, we have
\begin{equation}\label{add29}
\begin{aligned}
&\mathrm{LHS} \eqref{uk1}
\lesssim \|\frac{|\tau|^{\frac{1}{4}}}{|\tau|}\int_{|\xi-\eta|< |\tau|^{\frac{1}{2}}} |\mathcal{F}f| d(\xi-\eta)\|_{L^2(d \tau d(\xi+\eta))}\\
&\lesssim A \big\|{|\partial_t|^{\frac 1p-\frac12}|\nabla_y|^{\frac12}}{|\nabla_x|}^{\frac 12-\frac 2p}f\big\|_{L^2(dt)d(x-y)d(x+y)}
\end{aligned}\end{equation}
where
\begin{equation}\label{add30}
\begin{aligned}
A= \sup_{\tau, \xi+\eta}\frac{|\tau|^{\frac{3}{4}-\frac 1p}}{|\tau|} \left(\int_{|\xi-\eta|< |\tau|^{\frac{1}{2}}}
\frac{|\xi|^{\frac 4p-1}}{|\eta|}d(\xi-\eta)\right)^{\frac{1}{2}}.
\end{aligned}
\end{equation}
Changing variables, this is something like
\begin{align*}
A= \sup_{\tau, |u|<|\tau|^{\frac{1}{2}}}\frac{|\tau|^{\frac{3}{4}-\frac 1p}}{|\tau|} \left(\int_{|v|< |\tau|^{\frac{1}{2}}}
\frac{|u+v|^{\frac 4p-1}}{|u-v|}d v\right)^{\frac{1}{2}}.
\end{align*}
After a change of variables this is reduced to $\tau=1$, and $A$ is bounded.

Since by Sobolev, 
\begin{equation}\label{add31}
\begin{aligned}
 \big\||\partial_t|^{\frac 1p-\frac12} &|\nabla_y|^{\frac12}{|\nabla_x|}^{\frac12-\frac 2p}f\big\|_{L^2(dt)d(x-y)d(x+y)} \\
 =& \big\||\partial_t|^{\frac12-\frac 1p} {|\nabla_x|}^{\frac3q-\frac32}|\nabla_y|^{\frac12}{|\nabla_x|}^{\frac12}f\big\|_{L^2(dtdxdy)} 
 \\ \lesssim  & \| |\nabla_{x}|^{\frac12}|\nabla_{y}|^{\frac12} f\|_{L^{p'}(dt)L^{q'}(dx)L^2(dy) }.
\end{aligned}
\end{equation}
Thus the proof of \eqref{uk1} is complete, and the proof of \eqref{uk12} is similar.

\bigskip

Case 3:  $ |\xi|+ |\eta| > 2 |\tau|^{\frac{1}{2}}$ 

In this case, due to Case 1, we can assume additionally that $|\tau|^{\frac12}> 10(1+|\xi+\eta|)$.
Thus $|\xi-\eta|> |\xi+\eta|$, so also $|\xi-\eta|> |\tau|^{\frac{1}{2}}$. As before, it suffices to show \eqref{uk1} and \eqref{uk12}.
\begin{equation}\label{add32}
\begin{aligned}
&\mathrm{LHS} \eqref{uk1}
\lesssim \|\int_{2|\xi-\eta|> |\xi+\eta|+ |\tau|^{\frac{1}{2}}} |\frac{|\tau|^{\frac{1}{4}}|\mathcal{F}f|}{|\xi-\eta|^2} d(\xi-\eta)\|_{L^2(d \tau d(\xi+\eta))}\\
&\lesssim  A \big\|{|\partial_t|^{\frac 1p-\frac12}|\nabla_y|^{\frac12}}{|\nabla_x|}^{\frac 12-\frac 2p}f\big\|_{L^2(dt)d(x-y)d(x+y)}
\end{aligned}\end{equation}
In this case,
\begin{align}
A^2=\sup_{\xi+\eta, \tau}\int_{2|\xi-\eta|> |\xi+\eta|+ |\tau|^{\frac{1}{2}}} \frac{|\tau|^{\frac{3}{2}-\frac2p}}{|\xi-\eta|^4}\frac{|\xi|^{\frac4p-1}}{|\eta|} d(\xi-\eta)
\end{align}
Again we scale to $|\tau|^{\frac{1}{2}}+ |\xi+\eta|=1$ and have to estimate
\begin{align*}
\int_{|v|>1} \frac{1}{|v|^4}\frac{|u+v|^{\frac 4p-1}}{|u-v|} dv
\end{align*}
This is bounded uniformly in $|u|<1$. As before, the rest of the proof follow from Sobolev's inequality, and the proof of \eqref{uk12} is similar.

 \bigskip

Case 4 :$ \frac{1}{2}(|\xi|+|\eta|)<|\tau |^{\frac{1}{2}} < 2 (|\xi|+|\eta|)$, 
In this case, due to Case 1, we can assume additionally that $|\tau|^{\frac12}> 10(1+|\xi+\eta|)$
 so in this case $|\tau |^{\frac{1}{2}} \sim |\xi|\sim |\eta|\sim |\xi-\eta|$.
 
We shall use the decomposition $u=\sum_{k=0}^\infty P_{|\tau|\sim 2^{k}}u$, and  the square function estimate
\begin{equation}\label{add33}
\begin{aligned}
\|\dt u\|&_{L^{\infty}(d(x-y)) L^2(d(x+y)dt)} \\ 
&\sim  \|\Big(\sum_{k=0}^\infty |P_{|\tau|\sim 2^{k}} \dt u|^2\Big)^{\frac12}\|_{L^{\infty}(d(x-y)) L^2(d(x+y)dt)} \\
&\lesssim\Big(\sum_{k=0}^\infty 2^{2k} \|P_{|\tau|\sim 2^{k}}u\|_{{L^{\infty}(d(x-y)) L^2(d(x+y)dt)}}^2\Big)^{\frac12}.
\end{aligned}\end{equation}

For each fixed dyadic piece $P_{|\tau|\sim 2^{k}}u$, in the current case, we have 
$$P_{|\tau|\sim 2^{k}}u= P_{|\tau|\sim 2^{k}}P_{|\xi|\sim 2^{k/2}}P_{|\eta|\sim 2^{k/2}}u,$$
which implies 
 \begin{align*}
\|\big|P_{|\tau| \sim 2^k}u\|_{collapsing} &\lesssim
\|\dx P_{|\xi|\sim 2^{k/2}}P_{|\eta|\sim 2^{k/2}} u\|_{L^2(dt)L^6(d(x-y)) L^2(d(x+y))}\\ 
& \lesssim
2^{-k/2}\|\dx\dy P_{|\eta|\sim 2^{k/2}} u\|_{L^2(dt)L^6(d(x-y)) L^2(d(x+y))}
\end{align*}
where we used Bernstein's inequality in rotated coordinates twice, see e.g., Lemma 3.1 in \cite{1} for more details.

Thus,
\begin{equation}\label{add34}
 \begin{aligned}
&\|\dt u\|_{L^{\infty}(d(x-y)) L^2(d(x+y)dt)} \\
&\lesssim
\Big(\sum_{k=0}^\infty \|\dx\dy P_{|\eta|\sim 2^{k/2}} u\|_{L^2(dt)L^6(d(x-y)) L^2(d(x+y))}^2\Big)^{\frac12}
\\
&\lesssim
\left(\sum_{k=0}^\infty 
\| \dx\dy P_{|\eta| \sim 2^{k/2}}f\|^2_{\mathcal S'}\right)^{\frac{1}{2}}\\
& \lesssim \|\dx\dy f\|_{\mathcal S'}
\end{aligned}\end{equation}
where we used Strichartz (Theorem~\ref{2.1'}) in the second line, and square function estimates in $y$ the the last line.

\end{proof}

\newsection{Proof of Theorem~\ref{main}}

In this section, we shall see how we can apply the theorems in the previous section to prove Theorem~\ref{main}.
Throughout this section, besides \eqref{vcondition}, we shall use extensively the fact that
$$
\supp \hat{v}\subset B_1(0)$$
This assumption implies that, the multiplication operator $Tu=v_N u$ can at most enlarge the Fourier support of $u$ by a set of size $N$, which will greatly simplify our proof, especially the proof of Theorem~\ref{low collapsing1} below. 
If we assume $v$ satisfies \eqref{vcondition} without this condition,  one can follow similar steps in this section to get the same conclusion. In that case, the multiplication operator $Tu=v_N u$ can enlarge the Fourier support of $u$ by a set of arbitrary large size, but with a rapidly decay constant if the new Fourier support deviates from the Fourier support of $u$ by a large distance. For the sake of  simplicity, we do not present the full details here for this general case.

To begin with, we shall first prove the following theorem involving collapsing norms at low frequency

\begin{theorem}\label{low collapsing1}
Let $\Lambda$ satisfy \eqref{2.1},
we have
\begin{equation}
\begin{aligned}
& \|\dxx \Lm\|_{\lc}+ \|\dt \Lm\|_{\lc} \\
& \ls \|\dx\dy G\|_{\mathcal{S}'_{r}} 
+\e  \|\dxx H\|_{collapsing} +\e  \|\dt H\|_{collapsing}  \\
&\quad +\|\dx\dy \Lm_0\|_{L^2}
\end{aligned}
\end{equation}
where the norms $\|\cdot\|_{collapsing}$ and  $\|\cdot\|_{\lc}$ are defined as in \eqref{full collapsing} and \eqref{low collapsing}, and $\e$ is defined as in \eqref{vcondition}.
\end{theorem}
\begin{remark}\label{low open}
Due to the criticalness of the potential $N^2v(N(x-y))$ in \eqref{2.1}, it is still open to us if one can prove the above theorem without the frequency assumption on $\Lm$, i.e., to replace on the left side the $\|\cdot\|_{\lc}$ by the full collapsing norm. However, in the case where $\Lambda$ satisfy \eqref{2.1} with $H=0$, we do know how to control the full collapsing norm 
$$ \|\dxx \Lm\|_{collapsing}+ \|\dt \Lm\|_{collapsing}$$ by using a different argument, 
the details are given in Section 8. Also in \cite{1}, a stronger version of the theorem is proved for the case when one replace $N^2v(N(x-y))$ by $N^{3\beta-1}v(N^\beta(x-y))$ for some $\beta<1$.
\end{remark}
\begin{proof}
We shall focus on $\|\dt \Lm\|_{\lc}$, the proof for the term $\|\dxx \Lm\|_{\lc} $ is similar. And we will first show that 
\begin{equation}\label{dt collapsing}
\begin{aligned}
&\big\|P_{|\xi-\eta|<20N}\dt\Lambda\big\|_{collapsing} \\
& \ls \|\dx\dy G\|_{\mathcal{S}'_{r}} 
+\e  \|\dt H\|_{collapsing} +\|\dx\dy \Lm_0\|_{L^2} .
\end{aligned}
\end{equation}

Let $\rho\in C_0^\infty(\R^3)$ be a smooth partition of unity, which satisfies
$$ \sum_{j\in \mathbb{Z}^3} \rho(\xi-j)\equiv 1, \forall \,\xi\in \R^3.
$$
We also assume that $0\le \rho\le 1$,  $\rho\equiv 1$ if $|\xi|\le \frac12$, and $\supp \rho\in B_1(0)$, the unit ball in $\R^3$ centered at origin.  Let $\psi_j(x-y)$  be the inverse Fourier transform of $\rho(\frac{\xi-\eta}{40N}-j)$. For each fixed $j\in \mathbb{Z}^3$, 
 $\psi_j$ is Fourier supported in a ball of radius $40N$ centered at $40N\cdot j$. Denote
$$ P_{{j}} f=\psi_j*f
$$
so that the Fourier transform of $P_{j}f=\rho(\frac{\xi-\eta}{40N}-j) \hat f$. In particular,  let $P_0 f$ denote projection onto ball of radius $40N$ centered at origin, so that the Fourier transform of $P_{0}f=\rho(\frac{\xi-\eta}{40N}) \hat f$.

Define 
\begin{equation}\label{full norm collapsing}
\begin{aligned}
\|\dt\Lambda\|_{\mathcal{N}}=\sum_{k=0}^{N}2^{-k}\Big(\sum_{k\le |j|< k+1}\big\|P_j\dt\Lambda\big\|_{collapsing}\Big).
\end{aligned}
\end{equation}
It is clear that 
$$ \big\|P_{|\xi-\eta|<20N}\dt\Lambda\big\|_{collapsing} \ls  \big\|P_0\dt\Lambda\big\|_{collapsing}  \ls \|\dt\Lambda\|_{\mathcal{N}}.$$
Thus it suffices to show that 
\begin{equation}\label{redt collapsing}
\begin{aligned}
& \|\dt\Lambda\|_{\mathcal{N}} \\
& \ls \|\dx\dy G\|_{\mathcal{S}'_{r}} 
+\e  \|\dt H\|_{collapsing} +\|\dx\dy \Lm_0\|_{L^2}.
\end{aligned}
\end{equation}
We shall first deal with the last term in the norm, where $k=N$. 
Note that for each fixed $j\in \mathbb{Z}^3$, the Fourier transform in $x-y$ direction of $P_j\dt\Lambda$ is supported in a ball of radius $40N$, 
by Bernstein's  inequality, 
\begin{equation}\label{dt berstein}
\begin{aligned}
\|P_j \dt\Lambda\|&_{L^\infty (d(x-y))L^2(dt) L^2(d(x+y))}\\
\ls &N^{\frac12}\|P_j \dt\Lambda\|_{L^6 (d(x-y))L^2(dt) L^2(d(x+y))} \\
\ls &N^{\frac12}\| \dt\Lambda\|_{L^2(dt)L^6 (d(x-y)) L^2(d(x+y))}. \end{aligned}
\end{equation}

Recall that 
\begin{equation}\nonumber
\begin{aligned}
\S\Lambda(t,x,y)&=N^2v\big(N(x-y)\big)\Lambda(t,x,y)+G(t,x,y)\\
&\qquad+N^2v\big(N(x-y))H(t,x,y) \\
 \Lambda(0,\cdot)&=\Lambda_0.
\end{aligned} 
\end{equation}
By Theorem~\ref{Spsi2}, we have 
\begin{equation}\label{dt berstein1}
\begin{aligned}
\| \dt\Lambda\|&_{L^2(dt)L^6 (d(x-y)) L^2(d(x+y))}\\ 
&\le C \big\|N^2v\big(N(x-y)\big)\dt\Lambda \big\|_{L^2(dt)L^{\frac65}(d(x-y))L^2(d(x+y))}, \\
&\quad+C \big\|N^2v\big(N(x-y)\big)\dt H\big\|_{L^2(dt)L^{\frac65}(d(x-y))L^2(d(x+y))}, \\
&\qquad +C\big\|\big(\langle\nabla_{x-y}\rangle^{\frac{1}{2}}+\langle\nabla_{x+y}\rangle^{\frac{1}{2}}\big)G\big\|_{\sr} \\
&\qquad +C \big\|\big(\langle\nabla_{x-y}\rangle^{\frac{1}{2}}+\langle\nabla_{x+y}\rangle^{\frac{1}{2}}\big)\Lambda_0\big\|_{L^2}, \\
&\le C\e \big\|\dt\Lambda \big\|_{L^2(dt)L^{6}(d(x-y))L^2(d(x+y))}, \\
&\quad+C\e N^{-1/2} \big\|\dt H\big\|_{L^{\infty}(d(x-y))L^2(dt)L^2(d(x+y))}, \\
&\qquad +C\big\|\big(\langle\nabla_{x-y}\rangle^{\frac{1}{2}}+\langle\nabla_{x+y}\rangle^{\frac{1}{2}}\big)G\big\|_{\sr} \\
&\qquad +C \big\|\big(\langle\nabla_{x-y}\rangle^{\frac{1}{2}}+\langle\nabla_{x+y}\rangle^{\frac{1}{2}}\big)\Lambda_0\big\|_{L^2}, 
 \end{aligned} 
\end{equation}
where we used H\"older in the second inequality. By choosing $\e$ small enough such that $C\e<1/2$, we have 
\begin{equation}\label{dt berstein2}
\begin{aligned}
\| \dt\Lambda\|&_{L^2(dt)L^6 (d(x-y)) L^2(d(x+y))}\\ 
&\le  C\e N^{-1/2} \big\|\dt H\big\|_{L^{\infty}(d(x-y))L^2(dt)L^2(d(x+y))}, \\
&\qquad +C\big\|\big(\langle\nabla_{x-y}\rangle^{\frac{1}{2}}+\langle\nabla_{x+y}\rangle^{\frac{1}{2}}\big)G\big\|_{\sr} \\
&\qquad +C \big\|\big(\langle\nabla_{x-y}\rangle^{\frac{1}{2}}+\langle\nabla_{x+y}\rangle^{\frac{1}{2}}\big)u_0\big\|_{L^2}, \\
&\le  C\e N^{-1/2} \big\|\dt H\big\|_{L^{\infty}(d(x-y))L^2(dt)L^2(d(x+y))}, \\
&\qquad +C\big\|\dx \dy G\big\|_{\sr} +C \big\|\dx \dy u_0\big\|_{L^2}.
 \end{aligned} 
\end{equation}
If we combine \eqref{dt berstein} and \eqref{dt berstein2}, we get
\begin{equation}\label{dt berstein3}
\begin{aligned}
\|P_j \dt\Lambda\|&_{L^\infty (d(x-y))L^2(dt) L^2(d(x+y))}\\
&\le  C\e  \big\|\dt H\big\|_{L^{\infty}(d(x-y))L^2(dt)L^2(d(x+y))}, \\
& +CN^{\frac12}\big\|\dx \dy G\big\|_{\sr} +CN^{\frac12} \big\|\dx \dy u_0\big\|_{L^2}.
 \end{aligned}
\end{equation}

Since the number of $j\in \mathbb{Z}^3$ such that  $N\le|j|<N+1$ is bounded by $CN^2$,
we have 
\begin{equation}\label{dt berstein4}
\begin{aligned}
\sum_{N\le |j|< N+1}&2^{-N}\big\|P_j\dt\Lambda\big\|_{collapsing}\\
&\le  C2^{-N}N^2\e  \big\|\dt H\big\|_{L^{\infty}(d(x-y))L^2(dt)L^2(d(x+y))}, \\
&\quad +C2^{-N}N^{\frac52}\big\|\dx \dy G\big\|_{\sr} \\
&\quad\qquad+C2^{-N}N^{\frac52} \big\|\dx \dy u_0\big\|_{L^2}.
\\
&\le  C\e  \big\|\dt H\big\|_{L^{\infty}(d(x-y))L^2(dt)L^2(d(x+y))}, \\
&\quad +C\big\|\dx \dy G\big\|_{\sr} +C \big\|\dx \dy u_0\big\|_{L^2}.
 \end{aligned}
\end{equation}

Now we shall control the terms $0\le k\le N-1$ using a bootstrap argument.
Write $\Lambda=\Lambda_1+\Lambda_2+\Lambda_3$, where
\begin{align*}
 \S \,\Lambda_1&= G,\,\, \,\,\, \Lambda_1(0,\cdot)=\Lambda_0 \\
  \S \, \Lambda_2&= N^2v\big(N(x-y)\big)H(t,x,y), \,\,\,\,\text{with initial conditions 0} \\
 \S \, \Lambda_3&= N^2v\big(N(x-y)\big)\Lambda(t,x,y), \,\,\,\,\text{with initial conditions 0}.
\end{align*}
By Lemma~\ref{xco3}, we have
\begin{equation}\label{la3}
\begin{aligned}
\sum_{k=0}^{N-1}&2^{-k}\Big(\sum_{k\le |j|< k+1}\big\|P_j\dt\Lambda_1\big\|_{collapsing}\Big) \\
\ls& \sum_{k=0}^{N-1}2^{-k}\Big(\sum_{k\le |j|< k+1}\big\|\dt\Lambda_1\big\|_{collapsing}\Big) \\
\ls& \sum_{k=0}^{N-1}2^{-k}\sum_{k\le |j|< k+1} \Big(\|\dx\dy G\|_{\mathcal{S}^\prime_{r}}+ \|\dx\dy \Lambda_0\|_{L^2}\Big) \\
\ls& \|\dx\dy G\|_{\mathcal{S}^\prime_{r}}+ \|\dx\dy \Lambda_0\|_{L^2}.
\end{aligned}
\end{equation}
And if we repeat the argument in \eqref{dt berstein}-\eqref{dt berstein3}, for each fixed $j$, we have
\begin{equation}\label{dt berstein4'}
\begin{aligned}
\|P_j \dt\Lambda_2\|&_{L^\infty (d(x-y))L^2(dt) L^2(d(x+y))}\\ 
&\le  C\e  \big\|\dt H\big\|_{L^{\infty}(d(x-y))L^2(dt)L^2(d(x+y))}.
 \end{aligned}
\end{equation}
Thus, 
\begin{equation}\label{la2}
\begin{aligned}
\sum_{k=0}^{N-1}&2^{-k}\Big(\sum_{k\le |j|< k+1}\big\|P_j\dt\Lambda_2\big\|_{collapsing}\Big) \\
\ls& \sum_{k=0}^{N-1}\e2^{-k}\Big(\sum_{k\le |j|< k+1}\big\|\dt H\big\|_{collapsing}\Big) \\
\ls& \e \big\|\dt H\big\|_{collapsing}.
\end{aligned}
\end{equation}
It remains to control the terms involving $\Lambda_3$, note that since $\hat v$ is supported in the unit ball centered at origin, we have 
\begin{equation}\nonumber
\begin{aligned}
\S \, P_j\Lambda_3&= P_jN^2v\big(N(x-y)\big)\Lambda(t,x,y) \\
&= P_jN^2v\big(N(x-y)\big)\sum_{j':|j'-j|\le 1} P_{j'}\Lambda(t,x,y).
\end{aligned}
\end{equation}
Thus we have the following analog of \eqref{dt berstein4'}
\begin{equation}\label{dt berstein5}
\begin{aligned}
\|P_j \dt\Lambda_3\|&_{L^\infty (d(x-y))L^2(dt) L^2(d(x+y))}\\ 
&\le  C\sum_{j':|j'-j|\le 1} \e  \big\|P_{j'}\dt\Lambda\big\|_{L^{\infty}(d(x-y))L^2(dt)L^2(d(x+y))},
 \end{aligned}
\end{equation}
Thus,
\begin{equation}\label{la1}
\begin{aligned}
\sum_{k=0}^{N-1}&2^{-k}\Big(\sum_{k\le |j|< k+1}\big\|P_j\dt\Lambda_3\big\|_{collapsing}\Big) \\
\ls& \sum_{k=0}^{N-1}\e2^{-k}\Big(\sum_{k\le |j|< k+1}\sum_{j':|j'-j|\le 1}\big\|P_{j'}\dt \Lambda\big\|_{collapsing}\Big). 
\end{aligned}
\end{equation}
By choosing $\e$ small enough, the right hand side of \eqref{la1} is bounded by $\frac12\|\dt\Lambda\|_{\mathcal{N}}$, thus the proof of \eqref{redt collapsing} is complete.

To prove \eqref{dt collapsing} when $P_{|\xi-\eta|<20N}$ is replace by $P_{|\xi|<20N}$ or $P_{|\eta|<20N}$, we use a complete similar argument, the only necessary change is to replace the use of Bernstein's  
inequality in \eqref{dt berstein} by Bernstein's  
inequality in rotated coordinates.

\end{proof}
\noindent$Proof \,\,of\,\, Theorem~\ref{main}.$

To prove Theorem~\ref{main}, it remains to control $\|\Lm\|_{\mathcal S_{x,y}}$, where $\|\cdot\|_{\mathcal S_{x,y}}$ is defined as in \eqref{restr}.

\noindent\textbf{Case 1}. Let us first assume $|\xi+\eta|\ge\frac{N}{10}$, the norms we are going to control are different for the case $|\xi+\eta|< \frac{N}{10}$, but in both cases the norms contains $\|\Lm\|_{\mathcal S_{x,y}}$. In this case,
the norm with respect to which the potential is a perturbation should be
\begin{align}
\|\Lambda\|_{\mathcal N}&= \|\dx\dy\Lambda\|_{\mathcal S} \label{1655} \\ &
+\|P_{|\xi|\ge10N}P_{|\eta|<10N}\dx\Lambda\|_{collapsing}  \label{co11}\\
&+\|P_{|\xi|<10N}P_{|\eta|\ge10N}\dy\Lambda\|_{collapsing}\label{co22}
\end{align}
where $\|\cdot\|_{\mathcal{S}}$ is defined in \eqref{restrx-y} and  $\|\cdot\|_{collapsing}$ is defined as in \eqref{full collapsing}. Here we are abusing notations a bit by using $\Lambda$ to denote $P_{|\xi+\eta|\ge \frac{N}{10}}\Lambda$. The projection on $|\xi+\eta|\ge \frac{N}{10}$ is necessary, as we shall see later in the proof, we do not know how to control the collapsing norms in \eqref{co11} and \eqref{co22} without this assumption. 

Let's first control the low frequency, let $P_{<N}=P_{|\xi|<10N}P_{|\eta|<10N}$ and  $P_{>N}= I- P_{<N}$. Note that
\begin{equation}\label{4add1}
\begin{aligned}
\S \dx&\dy P_{<N}\Lambda\\
=&P_{<N}\dx\dy N^2v(N(x-y))P_{|\xi-\eta|<20N} \Lambda \\
& +\dx\dy P_{<N} G(t,x,y)\\
&+P_{<N}\dx\dy N^2v(N(x-y)) H(t,x,y)  \\
\sim &P_{<N} N^3v(N(x-y))P_{|\xi-\eta|<20N} \Lambda \\
& +\dx\dy P_{<N} G(t,x,y) \\
& +P_{<N} N^3v(N(x-y)) H(t,x,y)  
\end{aligned}
\end{equation}
Here we are abusing the notation a bit by writing
$$P_{<N}\dx\dy N^2v(N(x-y))P_{|\xi-\eta|<20N}  \sim P_{<N} N^3v(N(x-y))P_{|\xi-\eta|<20N},
$$
since, as a result of Bernstein's  inequality, for all $1\le p, q\le \infty$, we have
\begin{equation}\label{4add2}
\begin{aligned}
\| P_{<N}\dx\dy& N^2v(N(x-y))P_{|\xi-\eta|<20N} \Lambda\|_{L^p(x-y)L^q(x+y)}\\
\ls &\|P_{<N} N^3v(N(x-y))P_{|\xi-\eta|<20N} \Lambda \|_{L^p(x-y)L^q(x+y)},
\end{aligned}\end{equation}
which is harmless for our purposes. We shall use the same notation $\sim$ repeatly in the later arguments. Also, strictly speaking, 
$$P_{<N}N^2v(N(x-y))\Lambda=P_{<N} N^2v(N(x-y))P_{|\xi-\eta|<21N}\Lambda,
$$
instead of $P_{|\xi-\eta|<20N}\Lambda$, due to the fact that convolution with $\hat v$ will shift the frequency support. But it will not make a essential difference
in our argument.

By Strichartz estimate (Theorem~\ref{2.1'}), we have 
\bgd\label{420}\ad
&\| \dx\dy P_{<N}\Lambda \|_{\mathcal S} \\
&\ls   \|P_{|\xi-\eta|<20N} N^3v(N(x-y)) \Lambda\|_{L^2(dt)L^{\frac65}{(d(x-y))}L^2(d(x+y))}  \\
&\quad+  \|P_{|\xi-\eta|<20N} N^3v(N(x-y)) H\|_{L^2(dt)L^{\frac65}{(d(x-y))}L^2(d(x+y))}  \\
& \qquad + \| \dx\dy P_{<N}G \|_{\sr}+  \| \dx\dy \Lambda_0 \|_{L^2} \\
&\lesssim \|v\|_{L^{\frac65}}
\Big( \|P_{|\xi-\eta|<20N}\dxx\Lambda\|_{collapsing}+ \|\dxx H \|_{collapsing} \Big) 
\\ & \qquad + \| \dx\dy P_{<N}G \|_{\sr}+  \| \dx\dy \Lambda_0 \|_{L^2},
\ead\ed
where in the second inequality we used Bernstein's inequality and the fact that $|\xi+\eta|\ge \frac{N}{10}$.

By Theorem \ref{low collapsing1}, we have 
\begin{equation}\label{4add3}
\begin{aligned}
\| \dx\dy P_{<N}\Lambda \|_{\mathcal S} \lesssim\e &\|\dxx H \|_{collapsing} \\
 &+ \| \dx\dy P_{<N}G \|_{\mathcal{S}^\prime_{r}}+  \| \dx\dy \Lambda_0 \|_{L^2}. 
\end{aligned} \end{equation}

At high frequency, write
\begin{equation}\label{decom}
\begin{aligned}
 P_{>N}\Lambda=&P_{|\xi|\ge10N}P_{|\eta|\ge10N}\Lambda+P_{|\xi|\ge10N}P_{|\eta|< 10N}\Lambda \\
 &\hspace{2cm}+P_{|\xi|<10N}P_{|\eta|\ge10N}\Lambda\\
 =& I+II+III. 
 \end{aligned}
\end{equation}
\bigskip

To handle the first term $I$, note that
\begin{equation}\label{4add4}
\begin{aligned}
\S \dx&\dy P_{|\xi|\ge10N} P_{|\eta|\ge10N}\Lambda \\
 \sim& N^2v(N(x-y)) \dx\dy{P}_{|\xi|\ge9N} {P}_{|\eta|\ge 9N}  \Lambda  \\
 &+  \dx\dy P_{|\xi|\ge10N} P_{|\eta|\ge 10N}G(t,x,y)  \\
 &+ N^2v(N(x-y)) \dx\dy{P}_{|\xi|\ge9N} {P}_{|\eta|\ge 9N} H.
\end{aligned}
\end{equation}
The lower bounds on $|\xi|$ and $|\eta|$ changed slightly due to convolution with $\hat{v}_N$, since
$\hat{v}_N$ is compact supported in a set of size $N$.

Using Strichartz(Theorem~\ref{2.1'}) and H\"older's inequality, it is not hard to see that
\begin{equation}\label{4add5}
\begin{aligned}
\| &\dx\dy P_{|\xi|\ge10N} P_{|\eta|\ge10N}\Lambda \|_{\mathcal S} \\
&\lesssim\e \| \dx\dy \Lambda \|_{\mathcal S} \\
&\quad+\e \| \dx\dy H \|_{L^2(dt)L^6{(d(x-y))}L^2(d(x+y))} 
\\ & \qquad + \| \dx\dy G \|_{\sr}+  \| \dx\dy \Lambda_0 \|_{L^2},
\end{aligned} \end{equation}
where the $\e$ comes from $\|v\|_{L^{3/2}}$ when applying H\"older's inequality.
 
\bigskip

To handle the second term $II$, 
\begin{equation}\label{4add7}
\begin{aligned}
\S \dx&\dy P_{|\xi|\ge10N} P_{|\eta|<10N}\Lambda \\
 \sim &N^{\frac52}v(N(x-y)) \dx P_{|\xi|\ge10N} P_{|\eta|<10N} \Lambda  \\
& +N^{2}v(N(x-y)) \dx\dy P_{|\xi|\ge10N} P_{10N\le|\eta|<11N} \Lambda  \\
& + N^3v(N(x-y))  {P}_{9N\le|\xi|<10N} {P}_{|\eta|< 10N} \Lambda  \\
&+ N^{\frac52}v(N(x-y)) \dx{P}_{|\xi|\ge9N} {P}_{|\eta|<11N} H\\
 &+  \dx\dy P_{|\xi|\ge10N} P_{|\eta|< 10N}G(t,x,y).
 \end{aligned}\end{equation}
Again the bounds on $|\xi|$ and $|\eta|$ changed slightly due to convolution with $\hat{v}_N$.
By using Strichartz and H\"older's inequality, 
\begin{equation}\label{4add8}
\begin{aligned}
\| &\dx\dy P_{|\xi|\ge10N} P_{|\eta|<10N}\Lambda \|_{\mathcal S} \\
&\lesssim\e \| \dx P_{|\xi|\ge10N} P_{|\eta|<10N}\Lambda \|_{L^{\infty}{(d(x-y))}L^2(dt)L^2(d(x+y))} \\
&\quad+ \e \| \dx\dy \Lambda \|_{L^2(dt)L^{6}{(d(x-y))}L^2(d(x+y))} \\
&\quad+\| N^3v(N(x-y)) {P}_{9N\le|\xi|<10N} {P}_{ |\eta|< 10N} \Lambda  \|_{L^2(dt)L^{\frac65}{(d(x-y))}L^2(d(x+y))} \\
&\quad+\e \| \dx H \|_{L^{\infty}{(d(x-y))}L^2(dt)L^2(d(x+y))}
\\ & \quad + \| \dx\dy G \|_{\sr}+  \| \dx\dy \Lambda_0 \|_{L^2}.
\end{aligned} \end{equation}
Again, since we are assuming
$|\xi+\eta|\ge \frac{N}{10}$,  by Bernstein's inequality and H\"older's inequality, the third term on the right side can be controlled by 
\begin{align*}\nonumber
 \| &N^3v(N(x-y)) P_{9N\le |\xi|\le 10N}P_{ |\eta|\le 10N}  \Lambda\|_{L^2(dt)L^{6/5}{(d(x-y))}L^2(d(x+y))}  \\
 & \lesssim \|P_{|\xi-\eta|<20N}\dxx\Lambda\|_{collapsing}.\end{align*}
The third term $III$ can be handled in a similar way as the second term.

To finish to discussion for the case $|\xi+\eta|\ge \frac{N}{10}$,  we are reduced to estimate
\begin{equation}\label{4add9}
\begin{aligned}
&\|P_{|\xi|\ge10N}P_{|\eta|<10N}\dx\Lambda\|_{collapsing} \\
&+\|P_{|\xi|<10N}P_{|\eta|\ge10N}\dy\Lambda\|_{collapsing}.
\end{aligned}\end{equation}
We shall focus on the first term, 
since the other term involving $\dy$ can be dealt with similarly.

Write $\Lambda=\Lambda_1+\Lambda_2+\Lambda_3$, where
\begin{equation}\label{4add10}
\begin{aligned}
 \S \,\Lambda_1&= G,\,\, \,\,\, \Lambda_3(0,\cdot)=\Lambda_0\\
 \S \, \Lambda_2&= N^2v\big(N(x-y)\big)H(t,x,y), \,\,\,\,\text{with initial conditions 0} \\
  \S \, \Lambda_3&= N^2v\big(N(x-y)\big)\Lambda(t,x,y), \,\,\,\,\text{with initial conditions 0}.
\end{aligned}\end{equation}
By Lemma~\ref{xco1}, we have
\begin{equation}\label{la31}
\begin{aligned}
&\|P_{|\xi|\ge10N}P_{|\eta|<10N}\dx\Lambda_1\|_{collapsing}, \\
&\ls \|\dx\dy G\|_{\mathcal{S}^\prime_{r}}+ \|\dx\dy \Lambda_0\|_{L^2}.
\end{aligned}
\end{equation}

To handle $\Lambda_2$, since $|\eta|<10N$, by Bernstein's inequality at an angle, we have,
\begin{multline}
\|P_{|\xi|\ge10N}P_{|\eta|<10N}\dx \Lambda_2\|_{collapsing}\\ \lesssim N^{1/2} \|P_{|\xi|\ge10N}P_{|\eta|<10N}\dx\Lambda_2\|_{L^2(dt)L^{6}{(d(x-y))}L^2(d(x+y))}.
\end{multline}
Note that
\begin{equation}\label{4add11}
\begin{aligned}
\S  \dx & N^{1/2} P_{|\xi|\ge 10N}P_{|\eta|<10N}\Lambda_2 \\
& \sim N^{5/2}v(N(x-y)) \dx{P}_{|\xi|\ge 9N} {P}_{|\eta|<11N}  H.
\end{aligned}\end{equation}
By Strichartz and H\"older's inequality
\begin{equation}\label{4add12}
 \begin{aligned}
N&^{1/2} \|P_{|\xi|\ge10N}P_{|\eta|<10N}\dx\Lambda_2\|_{L^2(dt)L^{6}{(d(x-y))}L^2(d(x+y))}\\
&\ls \| N^{\frac52}v(N(x-y)) {P}_{|\xi|\ge 9N} {P}_{ |\eta|< 11N} H  \|_{L^2(dt)L^{\frac65}{(d(x-y))}L^2(d(x+y))} \\
&\lesssim\e \| \dx H \|_{L^{\infty}{(d(x-y))}L^2(dt)L^2(d(x+y))}.
\end{aligned} \end{equation}

To handle $\Lambda_3$, since $|\eta|<10N$, by Bernstein's inequality at an angle, we still have,
\begin{multline}
\|P_{|\xi|\ge10N}P_{|\eta|<10N}\dx \Lambda_3\|_{collapsing}\\ \lesssim N^{1/2} \|P_{|\xi|\ge10N}P_{|\eta|<10N}\dx\Lambda_3\|_{L^2(dt)L^{6}{(d(x-y))}L^2(d(x+y))}.
\end{multline}
Note that
\begin{equation}\label{4add13}
\begin{aligned}
\S  \dx & N^{1/2} P_{|\xi|\ge 10N}P_{|\eta|<10N}\Lambda_3 \\
& \sim N^{\frac52}v(N(x-y)) \dx  P_{|\xi|\ge 10N}P_{|\eta|<10N}  \Lm \\
& +N^{\frac52}v(N(x-y)) \dx P_{|\xi|\ge10N} P_{10N\le|\eta|<11N} \Lambda  \\
&+  N^3v(N(x-y)) {P}_{9N\le|\xi|<10N} {P}_{|\eta|< 10N} \Lambda. 
\end{aligned}\ed
By Strichartz and H\"older's inequality
\begin{equation}\label{4add14}
\begin{aligned}
N&^{\frac12}\| \dx P_{|\xi|\ge10N} P_{|\eta|<10N}\Lambda_3 \|_{L^2(dt)L^{6}{(d(x-y))}L^2(d(x+y))} \\
&\lesssim\e \| \dx P_{|\xi|\ge10N} P_{|\eta|<10N}\Lambda \|_{L^{\infty}{(d(x-y))}L^2(dt)L^2(d(x+y))} \\
&\quad+ \e N^{\frac12}\| \dx \Lambda \|_{L^2(dt)L^{6}{(d(x-y))}L^2(d(x+y))}\\
&\quad+\| N^3v(N(x-y)) {P}_{9N\le|\xi|<10N} {P}_{ |\eta|< 10N} \Lambda  \|_{L^2(dt)L^{\frac65}{(d(x-y))}L^2(d(x+y))},\\
&\lesssim\e \| \dx P_{|\xi|\ge10N} P_{|\eta|<10N}\Lambda \|_{L^{\infty}{(d(x-y))}L^2(dt)L^2(d(x+y))} \\
&\quad+ \e \| \dx\dy \Lambda \|_{L^2(dt)L^{6}{(d(x-y))}L^2(d(x+y))}\\
&\quad+\| {P}_{ |\xi-\eta|< 20N} \dxx\Lambda  \|_{collapsing},
\end{aligned} \ed
where for the second term on the right side we used Bernstein's inequality and the fact that $|\eta|\ge{10N}$
, and for the third term in the right side, we used Bernstein's inequality and the fact that $|\xi+\eta|\ge \frac{N}{10}$.

\noindent\textbf{Case 2}.  $|\xi+\eta|< \frac{N}{10}$,  $|\xi-\eta|< 10N$. 

In this case, we are only able to control $\|\Lm\|_{\mathcal{S}_{x,y}}$ instead of $\|\Lm\|_{\mathcal{S}}$, recall that $\|\cdot\|_{\mathcal{S}_{x,y}}$ is defined as in \eqref{restr}.
Note that
\begin{equation}\label{436}
\begin{aligned}
\S \dx&\dy P_{|\xi+\eta|< \frac{N}{10}}P_{|\xi-\eta|< 10N}\Lambda\\
\sim & P_{|\xi+\eta|< \frac{N}{10}}P_{|\xi-\eta|< 10N} N^3v(N(x-y))P_{|\xi-\eta|<20N} \Lambda \\
& +\dx\dy P_{|\xi+\eta|< \frac{N}{10}}P_{|\xi-\eta|< 10N}G(t,x,y) \\
& + P_{|\xi+\eta|< \frac{N}{10}}P_{|\xi-\eta|< 10N}N^3v(N(x-y)) H(t,x,y)  
\end{aligned}
\end{equation}

By using Theorem~\ref{Spsi} for the first and third term, and the Strichartz estimate (Theorem~\ref{2.1'}) for the second term on the right side of \eqref{436}, we have 
\bgd\ad
&\| \dx\dy P_{|\xi+\eta|< \frac{N}{10}}P_{|\xi-\eta|< 10N}\Lambda \|_{\mathcal{S}_{x,y}} \\
&\ls   \|P_{|\xi-\eta|<20N} N^3v(N(x-y))\dxx\Lambda\|_{L^{1}{(d(x-y))}L^2(dt)L^2(d(x+y))}  \\
&\quad+   \|P_{|\xi-\eta|<20N} N^3v(N(x-y)) \dt\Lambda\|_{L^{1}{(d(x-y))}L^2(dt)L^2(d(x+y))}  \\
&\quad+  \|P_{|\xi-\eta|<20N} N^3v(N(x-y)) \dxx H\|_{L^{1}{(d(x-y))}L^2(dt)L^2(d(x+y))}  \\
&\quad+  \|P_{|\xi-\eta|<20N} N^3v(N(x-y)) \dt H\|_{L^{1}{(d(x-y))}L^2(dt)L^2(d(x+y))}  \\
& \qquad + \| \dx\dy G \|_{\sr}+  \| \dx\dy \Lambda_0 \|_{L^2} \\
&\lesssim \|v\|_{L^{1}}
\Big( \|P_{|\xi-\eta|<20N}\dxx\Lambda\|_{collapsing}+\|P_{|\xi-\eta|<20N}\dt\Lambda\|_{collapsing} \\ 
&\quad\qquad\qquad+ \|\dxx H \|_{collapsing}+\|\dt H \|_{collapsing} \Big) 
\\ & \qquad + \| \dx\dy G \|_{\sr}+  \| \dx\dy \Lambda_0 \|_{L^2},
\ead\ed

By Theorem \ref{low collapsing1}, we have 
\begin{equation}
\begin{aligned}
\| &\dx\dy P_{|\xi+\eta|< \frac{N}{10}}P_{|\xi-\eta|< 10N}\Lambda \|_{\mathcal S} \lesssim\e \|\dxx H \|_{collapsing} \\
 &+ \e \|\dt H \|_{collapsing}+ \| \dx\dy G \|_{\mathcal{S}^\prime_{r}}+  \| \dx\dy \Lambda_0 \|_{L^2}. 
\end{aligned} \end{equation}

\vspace{.5cm}
\noindent\textbf{Case 3}.  $|\xi+\eta|< \frac{N}{10}$,  $|\xi-\eta|\ge 10N$. 

In this case, $|\xi|\sim|\eta|\sim |\xi-\eta|$, define $P_{>N}=P_{|\xi+\eta|< \frac{N}{10}} P_{|\xi-\eta|\ge 10N}$
The norm with respect to which the potential is a perturbation should be
\begin{align}
\|\Lambda\|_{\mathcal N}=& \||\dx\dy P_{>N}\Lambda\|_{\mathcal S_{x, y}} \\ &+
\||\dt\dyy\widetilde{P}_{>N}\Lambda\|_{L^2(dt)L^{6}(d(x-y))L^2(d(x+y))}\label{867} \\
&+\|\dxx\dyy\widetilde{P}_{>N}\Lambda\|_{L^2(dt)L^{6}(d(x-y))L^2(d(x+y))},
\end{align}
where $\widetilde{P}_{>N}=P_{|\xi+\eta|< \frac{N}{10}} P_{|\xi-\eta|\ge 9N}$.

In this case,
\begin{equation}\label{4add16}
\begin{aligned}
\S &\dx\dy P_{>N}\Lambda \\
& \sim N^2v(N(x-y))\langle\nabla_{x-y}\rangle\widetilde{P}_{>N} \Lambda \\
& +\dx\dy P_{>N} G(t,x,y) +N^2v(N(x-y)\langle\nabla_{x-y}\rangle\widetilde{P}_{>N} H.
\end{aligned}\ed
Using Strichartz 
 for the second and third terms on the right side, Theorem \ref{Spsi1} for the first term on the RHS, and H\"older's inequality, we get
 \begin{equation}\label{4add17}
\begin{aligned}
&\| \dx\dy P_{>N}\Lambda \|_{\mathcal S_{x, y}} \\
&\lesssim\e
\|\dt\dyy\widetilde{P}_{>N}\Lambda\|_{L^2(dt)L^{6}(d(x-y))L^2(d(x+y))}
\\ & \quad +\e \|\dxx\dyy\widetilde{P}_{>N}\Lambda\|_{L^2(dt)L^{6}(d(x-y))L^2(d(x+y))} \\
&\quad+ \| \dx\dy P_{>N}G \|_{\sr} \\
&\quad+\e\|\langle\nabla_{x-y}\rangle\widetilde{P}_{>N}H\|_{L^2(dt)L^{6}(d(x-y))L^2(d(x+y))} \\
&\lesssim\e
\|\dt\dyy\widetilde{P}_{>N}\Lambda\|_{L^2(dt)L^{6}(d(x-y))L^2(d(x+y))}
\\ & \quad +\e \|\dxx\dyy\widetilde{P}_{>N}\Lambda\|_{L^2(dt)L^{6}(d(x-y))L^2(d(x+y))} \\
&\quad+ \| \dx\dy P_{>N}G \|_{\sr} \\
&\quad+\e\|\dx\dy H\|_{L^2(dt)L^{6}(d(x-y))L^2(d(x+y))}. 
\end{aligned}\ed
Thus, it suffice to control 
 \begin{equation}
\begin{aligned}
\|&\dt\dyy\widetilde{P}_{>N}\Lambda\|_{L^2(dt)L^{6}(d(x-y))L^2(d(x+y))} \\
&+\||\dxx\dyy\widetilde{P}_{>N}\Lambda\|_{L^2(dt)L^{6}(d(x-y))L^2(d(x+y))} \\
=&I+II.
\end{aligned}\ed
For simplicity, we shall only give the details for the first term $I$, the second term is easier and can be handled in a similar way.

Note that
 \begin{equation}
\begin{aligned}
\S\dyy &P_{|\xi-\eta|\ge 9N}P_{|\xi+\eta|< \frac{N}{10}}\Lambda \\
 \sim &N^2v(N(x-y))\dyy{P}_{|\xi-\eta|\ge 9N} P_{|\xi+\eta|< \frac{N}{10}} \Lambda \\
& + N^{\frac52}v(N(x-y)){P}_{8N\le|\xi-\eta|<9N} P_{|\xi+\eta|< \frac{N}{10}} \Lambda \\
 &+\dyy P_{|\xi-\eta|\ge 9N}P_{|\xi+\eta|< \frac{N}{10}}G(t,x,y) \\
 &+ N^2v(N(x-y))\dyy{P}_{|\xi-\eta|\ge 8N} P_{|\xi+\eta|< \frac{N}{10}}H,
\end{aligned}\ed
where the bounds on $|\xi-\eta|$ changed slightly due to convolution with $\hat{v}_N$.

By Theorem~\ref{Spsi2} and H\"older's inequality, we have 
 \begin{equation}
\begin{aligned}
&\|\dt\dyy P_{|\xi-\eta|\ge 9N}P_{|\xi+\eta|< \frac{N}{10}}\Lambda\|_{L^2(dt)L^{6}(d(x-y))L^2(d(x+y))}\\
&\ls \e \|\dt\dyy{P}_{|\xi-\eta|\ge 9N} P_{|\xi+\eta|< \frac{N}{10}}\Lambda\|_{L^2(dt)L^{6}(d(x-y))L^2(d(x+y))} \\
&\,+  \||N^{\frac52}v(N(x-y)) {P}_{8N\le |\xi-\eta|<9N} P_{|\xi+\eta|< \frac{N}{10}}\dt \Lambda\|_{L^2(dt)L^{6}(d(x-y))L^2(d(x+y))} \\
&\, +\|\langle\nabla_{x-y}\rangle{P}_{|\xi-\eta|\ge 9N} P_{|\xi+\eta|< \frac{N}{10}}G\|_{\sr}\\
&\, +\e\|\langle\nabla_{x-y}\rangle\widetilde{P}_{|\xi-\eta|\ge 8N} P_{|\xi+\eta|< \frac{N}{10}}H\|_{L^2(dt)L^{6}(d(x-y))L^2(d(x+y))} \\
&\,+\|\langle\nabla_{x-y}\rangle\widetilde{P}_{|\xi-\eta|\ge 9N} P_{|\xi+\eta|< \frac{N}{10}}\Lambda_0\|_{L^2}\\
&\ls \e \|\dt\dyy\widetilde{P}_{>N} \Lambda\|_{L^2(dt)L^{6}(d(x-y))L^2(d(x+y))} \\
&\,+\e\|\dt{P}_{|\xi-\eta|\le20N}\Lambda\|_{L^{\infty}(d(x-y))L^2(dt)L^2(d(x+y))}\\
&\, +\|\dx\dy G\|_{\mathcal{S}^\prime_{r}}+\e\|\dx\dy H\|_{L^2(dt)L^{6}(d(x-y))L^2(d(x+y))}  \\
& \, \qquad+\|\dx\dy \Lambda_0\|_{L^2},
\end{aligned}\ed
where the first two terms on the right side of the first inequality corresponds to the first term on the right side of \eqref{dtdx}, and the last three terms on the right side of the first inequality corresponds to the remaining terms on the right side of \eqref{dtdx}.

\newsection{Estimates for the nonlinear equation}

Recall the notation
$$\S_{\pm}=\frac1i\frac{\partial}{\partial t}-\Delta_x+\Delta_y
$$
From now on, $V_N(x)=N^3v(Nx)$.

Define $\Gamma=\Gamma_c+\Gamma_p$, $\Lambda=\lcc+\lpp$, where $\gc=\bar{\phi}\otimes\phi$, $\lcc={\phi}\otimes\phi$, $\gp=\frac1N\overline{\sh}\circ \sh$,
and $\lpp=\frac{1}{2N}\mbox{sh} (2k)$. Let $\rho(t,x)=\Gamma(t,x,x)$.

The four relevant equations are
\begin{align}
\S&\lpp+\{V_N*\rho, \lpp\}+\frac{V_N}{N}\lpp \label{lambdap}\\
&+\big((V_N\bar\Gamma_p)\circ\lpp+(V_N\lpp)\circ\gp\big)_{symm} \nonumber\\
&+\big((V_N\bar\Gamma_c)\circ\lpp+(V_N\lcc)\circ\gp\big)_{symm}=-\frac{V_N}{N}\lcc \nonumber\\
\S&_{\pm}\gp+[V_N*\rho, \gp]+\big((V_N\gp)\circ\gp+(V_N\bar\Lambda_p)\circ\lpp\big)_{skew} \label{gammap} \\
&+\big((V_N\gc)\circ\gp+(V_N\bar\Lambda_c)\circ\lpp\big)_{skew}=0\nonumber \\
\S&\lcc+\{V_N*\rho, \lcc\}+\big((V_N\bar\Gamma_p)\circ\lcc+(V_N\lpp)\circ\gc\big)_{symm}=0 \label{lambdac}\\
\S&_\pm\gc+[V_N*\rho, \gc]+\big((V_N\Gamma_p)\circ\gc+(V_N\bar\Lambda_p)\circ\lcc\big)_{skew}=0 \label{gammac}.
\end{align}
Here $\big(A(x,y)\big)_{symm}=A(x,y)+A(y,x)$ and $\big(A(x,y)\big)_{skew}=A(x,y)-\bar A(y,x)$

The norm used for $\lpp$ is called $\mathcal{N}_1(\Lambda)$ and is 
\begin{equation}\label{lambdap norm}
\begin{aligned}
\|\Lm\|&_{\mathcal{N}_1(\Lambda)}=\|\dx\dy \Lm\|_{\mathcal{S}_{x, y}} + \|\dxx \Lm\|_{\lc}\\
 & + \|\dt \Lm\|_{\lc}\\
 &+\big\|\dx\dy P_{|\xi|\ge10N}P_{|\eta|\ge{10N}}P_{|\xi+\eta|\ge \frac{N}{10}}\Lm\big\|_{L^{2}(dt)L^6(d(x-y))L^{2}(d(x+y))}\\
 &+\big\|\dxx\dyy P_{|\xi-\eta|\ge10N}P_{|\xi+\eta|< \frac{N}{10}}\Lm\big\|_{L^{2}(dt)L^6(d(x-y))L^{2}(d(x+y))} \\
  &+\big\|\dt\dyy P_{|\xi-\eta|\ge10N}P_{|\xi+\eta|< \frac{N}{10}}\Lm\big\|_{L^{2}(dt)L^6(d(x-y))L^{2}(d(x+y))},
\end{aligned}
\end{equation}
where the norm $\|\cdot\|_{\lc}$ is defined as in \eqref{low collapsing}.  The last three norms in \eqref{lambdap norm} does not appear in the statement of Theorem~\ref{main}, but as 
one can see from the proof of Theorem~\ref{main} in the previous section, they satisfy the same bounds as the first three norms on the right side of \eqref{lambdap norm}.

The norm used for $\lcc$ is called $\mathcal{N}_2(\Lambda)$ and is 
\begin{equation}\label{lambdac norm}
\begin{aligned}
\|\Lm\|&_{\mathcal{N}_2(\Lambda)}=\|\dx\dy \Lm\|_{\mathcal{S}_{x, y}} + \|\dxx \Lm\|_{collapsing}\\
 & + \|\dt \Lm\|_{collapsing} + \|\dx \Lm\|_{collapsing}\\
 & + \|\dy \Lm\|_{collapsing} + \big\|\dx\dy \Lm\big\|_{L^{2}(dt)L^6(d(x-y))L^{2}(d(x+y))}.
\end{aligned}
\end{equation}

We will use the following a priori estimates for $\Gamma(t,x,x)$ (proved in Lemma 6.2 in \cite{1}).

\begin{lemma}\label{a priori}
Let the potential $v$ satisfies \eqref{vconditio}, and
the initial conditions satisfy \eqref{initial condition}, we have for all $0\le\alpha\le1$, 
\begin{equation}
\|\langle \nabla_{x+y}\rangle ^{\alpha}\Gamma\|_{L^8(dt)L^{\infty}(d(x-y))L^{\frac43}(d(x+y))} \ls 1 .
\end{equation}
The above estimates also hold for $\gp$, $\gc$ and $\lcc$ separately.
\end{lemma}

We need to use a continuity argument, we have to localize our estimates to intervals $[0,T]$, where the right end of the interval must be a variable $T$.
Define $\lcct, \lppt, \gpt, \gct$ to be solutions to the standard equations with the RHS multiplied by $\chi_{[0,T]}$:
\begin{align}
\S&\lppt+\frac{V_N}{N}\lppt \label{lambdap'}\\
&=\chi_{[0,T]}\Big(-\{V_N*\rho, \lpp\}-\big((V_N\bar\Gamma_p)\circ\lpp+(V_N\lpp)\circ\gp\big)_{symm} \nonumber\\
&\,\,\,-\big((V_N\bar\Gamma_c)\circ\lpp+(V_N\lcc)\circ\gp\big)_{symm}-\frac{V_N}{N}\lcc\Big) \nonumber\\
\S&_{\pm}\gpt=\chi_{[0,T]}\Big(-[V_N*\rho, \gp]-\big((V_N\gp)\circ\gp+(V_N\bar\Lambda_p)\circ\lpp\big)_{skew} \label{gammap'} \\
&-\big((V_N\gc)\circ\gp+(V_N\bar\Lambda_c)\circ\lpp\big)_{skew}\Big)\nonumber \\
\S&\lcct=\chi_{[0,T]}\Big(-\{V_N*\rho, \lcc\}-\big((V_N\bar\Gamma_p)\circ\lcc+(V_N\lpp)\circ\gc\big)_{symm} \Big)\label{lambdac'}\\
\S&_\pm\gct=\chi_{[0,T]}\Big(-[V_N*\rho, \gc]-\big((V_N\Gamma_p)\circ\gc+(V_N\bar\Lambda_p)\circ\lcc\big)_{skew}\Big) \label{gammac'}
\end{align}
with $\lcct(0,\cdot)=\lcc(0,\cdot)$,  and similarly for the other three functions. Also, we $\lcct=\lcc$ in $[0,T]$ (but not outside this interval), and similarly for the other three functions.

\begin{theorem}\label{lpb}
Let $[0,T]$ be as above, there exist a universal constant $C$ such that
\begin{equation}
\begin{aligned}
\|\lppt\|&_{\mathcal{N}_1(\Lambda)}\le C\|\dx\dy\Lambda_p(0, \cdot)\|_{L^2}+ C\e \|\lppt\|_{\mathcal{N}_1(\Lambda)} \\
&+C\e \|\gpt\|_{\mathcal{S}_{x,y}}+C\e \|\lppt\|_{\mathcal{N}_1(\Lambda)}\|\gpt\|_{\mathcal{S}_{x,y}}+C\e\|\lcct\|_{\mathcal{N}_2(\Lambda)}.
\end{aligned} 
\end{equation}
\end{theorem}
The proof is based on Theorem~\ref{main}, there exists a constant $C$ such that
\begin{align}
\|&\lppt\|_{\mathcal{N}_1(\Lambda)}\nonumber\\
& \le C\Big( \|\dx\dy\chi_{[0,T]}\big(\{V_N*\rho, \lpp\}+\big((V_N\bar\Gamma_p)\circ\lpp+(V_N\lpp)\circ\gp\big)_{symm}\nonumber \\
&\quad+\big((V_N\bar\Gamma_c)\circ\lpp+(V_N\lcc)\circ\gp\big)_{symm}\big)  \|_{\mathcal{S}'_{r}} \nonumber\\
&\quad+\e \|\dx\dy \chi_{[0,T]}\lcc\|_{L^2(dt)L^6(x-y) L^2(d(x+y))}\nonumber \\
&\quad+ \e  \|\dxx \chi_{[0,T]}\lcc\|_{collapsing} +\e  \|\dt \chi_{[0,T]}\lcc\|_{collapsing}\\
&\quad + \e  \|\dx \chi_{[0,T]}\lcc\|_{collapsing} +\e  \|\dy \chi_{[0,T]}\lcc\|_{collapsing}\nonumber\\
&\quad+ C\|\dx\dy\Lambda_p(0, \cdot)\|_{L^2}\Big).\nonumber
\end{align} 
For all terms other than $\|\dt \chi_{[0,T]}\lcc\|_{collapsing}$, the subscript $T$ can be trivially added to $\Lm$, $\Gamma$
on the RHS.  And we also have 
 \begin{equation}
\begin{aligned}
\|&\dt \chi_{[0,T]}\lcc\|_{L^{\infty}(d(x-y))L^2(dt)L^2(d(x+y))}\\
&=\|\dt \chi_{[0,T]}\lcct\|_{L^{\infty}(d(x-y))L^2(dt)L^2(d(x+y))}\\
&\ls  \|\dt \lcct\|_{L^{\infty}(d(x-y))L^2(dt)L^2(d(x+y))} \\
&\ls  \|\lcct\|_{\mathcal{N}_2(\Lambda)},
\end{aligned}\ed
where in the third line we used the fact that 
$$\|\dt \chi_{[0,T]}F\|_{L^{\infty}(d(x-y))L^2(dt)L^2(d(x+y))}\\
\ls  \|\dt F\|_{L^{\infty}(d(x-y))L^2(dt)L^2(d(x+y))},
$$
for any interval $[0,T]$. As remarked in \cite{1}, this can be shown by using the equivalent definition 
 \begin{equation}
\|\dt u\|_{L^2}=\int\int \frac{|u(t)-u(s)|^2}{|t-s|^{1+\frac12}} dtds,
\ed
and the generalized Hardy's inequality from \cite{4}.

In the lemma that follow, we estimate the norm of the nonlinear terms in suitable dual Strichartz norms, using the bound \eqref{a priori} whenever possible.
\begin{lemma}\label{nonlinear}
Let $[0,T]$ be as above, there exist a universal constant $C$ such that
\begin{equation}
\begin{aligned}
&\|\dx\dy\Big(\{V_N*\rho, \lppt\}+(V_N\bar\Gamma)\circ\lppt\Big)\|_{L^{\frac85}([0,T])L^{\frac43}(dx)L^2(dy)}. \\
&\le C\e \|\dx\dy\lppt\|_{\mathcal{S}_{x,y}}.\end{aligned} 
\end{equation}
Here $\Gamma$ can be $\gp$ or $\gc$. The result depends on the a priori bounds for $\Gamma$, but is true with $\lppt$ replaced with any other function.
\end{lemma}
\begin{proof}
In this case, we essentially view $V_N$ as a $\delta$ distribution, by Minkowski integral inequality, and it suffices to show that 
 \begin{equation}
\begin{aligned}
&\sup_z\|\dx\dy\Big(\Gamma(t,x,x+z)\lppt(t,x+z,y)\Big)\|_{L^{\frac85}([0,T])L^{\frac43}(dx)L^2(dy)} \\
&\le C \|\dx\dy\lppt\|_{\mathcal{S}_{x,y}},\end{aligned}\ed 
where the extra $\e$ factor in \eqref{nonlinear} can be remedied by the smallness of $\|V_N\|_{L^1}$.

Using the fractional Leiniz rule from Theorem 5.1 in \cite{1}, we have the following estimate, uniformly in $z$:
 \begin{equation}
\begin{aligned}
&\|\dx\dy\Big(\Gamma(t,x,x+z)\lppt(t,x+z,y)\Big)\|_{L^{\frac85}([0,T])L^{\frac43}(dx)L^2(dy)} \\
&\le C\|\dx\Gamma(t,x,x+z)\|_{L^{8}(dt)L^{\frac43+}(dx)} \|\dy\lppt\|_{L^{2}(dt)L^{\infty-}(dx)L^2(dy)} \\
&+C\|\Gamma(t,x,x+z)\|_{L^{8}(dt)L^{\frac{12}{7}}(dx)} \|\dx\dy\lppt\|_{L^{2}(dt)L^{6}(dx)L^2(dy)} \\
&\le C\|\langle \nabla_x\rangle^\alpha\Gamma(t,x,x+z)\|_{L^{8}(dt)L^{\frac43}(dx)} \|\dx\dy\lppt\|_{L^{2}(dt)L^{6}(dx)L^2(dy)} \\
&\le C  \|\dx\dy\lppt(t,x+z,y)\|_{L^{2}(dt)L^{6}(dx)L^2(dy)}.
\end{aligned}\ed
Here $\alpha$ is can be any number in $(\frac12,1]$, $\frac43+$ is a number that is bigger than but can be arbitrary close to $\frac43$, similarly $\infty-$ is any finite number but can
be arbitrary large. We use $\infty-$ since we do not have the sharp Sobolev  estimate from $L^6$ to $L^\infty$. In the last inequality we used Lemma~\ref{a priori}.
\end{proof}

Since $\lcc$ satisfies the same a priori estimates as $\Gamma$, by the exact same argument we get 
\begin{lemma}\label{nonlinear1}
Let $[0,T]$ be as above, there exist a universal constant $C$ such that
\begin{equation}
\begin{aligned}
&\|\dx\dy\Big((V_N\Lm_c)\circ\gpt\Big)\|_{L^{\frac43}([0,T])L^{\frac32}(dx)L^2(dy)} \\
&\le C\e  \|\dx\dy\gpt\|_{\mathcal{S}_{x,y}}.\end{aligned} 
\end{equation}
The result depends on the a priori bounds for $\lcc$, but is true with $\gpt$ replaced with any other function.
\end{lemma}

We continue estimating nonlinear terms.
\begin{lemma}\label{nonlinear2}
Let $[0,T]$ be as above, there exist a universal constant $C$ such that
\begin{equation}
\begin{aligned}
&\|\dx\dy\Big((V_N\lppt)\circ\gpt\Big)\|_{L^{\frac43}([0,T])L^{\frac32}(dx)L^2(dy)} \\
&\le C\e\|\lppt\|_{\mathcal{N}_1(\Lambda)} \|\dx\dy\gpt\|_{\mathcal{S}_{x,y}}.\end{aligned} 
\end{equation}
The result is still true if we replace $\gpt$ with any other function.
\end{lemma}
\begin{proof}
In this case, we shall not treat $V_N$ as a $\delta$ distribution.
Recall that
 \begin{equation} 
\begin{aligned}
\dx&\dy\Big((V_N\lppt)\circ\gpt\Big)\\
&=\dx\dy\int N^3v(N(x-z))\lppt(x,z)\gpt(z,y) dz \\
&=\int N^3v(Nz)\dx\dy\Big(\lppt(x,x+z)\gpt(x+z,y)\Big) dz.
\end{aligned} \ed
For fixed $z$, the following holds, uniformly in $z$
\begin{align*}
&\|\dx\dy\Big(\lppt(t,x,x+z)\gpt(t,x+z,y)\Big)\|_{L^{\frac43}([0,T])L^{\frac32}(dx)L^2(dy)} \\
&\le C\|\dx\lppt(t,x,x+z)\|_{L^{2}(dt)L^{2}(dx)} \|\dy\gpt(t,x,y)\|_{L^{4}(dt)L^{6}(dx)L^2(dy)} \\
&+C\|\lppt(t,x,x+z)\|_{L^{2}(dt)L^{3}(dx)} \|\dx\dy\gpt(t,x,y)\|_{L^{4}(dt)L^{3}(dx)L^2(dy)} \\
&\le C\|\dx\lppt(t,x,x+z)\|_{L^{2}(dt)L^{2}(dx)} \|\dx\dy\gpt(t,x,y)\|_{L^{4}(dt)L^{3}(dx)L^2(dy)}.
\end{align*}
Thus, by Minkowski integral inequality it suffices to show that 
\begin{align}
 \int |N^3v(Nz)| \|\dx\lppt(t,x,x+z)\|_{L^{2}(dt)L^{2}(dx)} dz
 \le  C\e\|\lppt\|_{\mathcal{N}_1(\Lambda)},
\end{align}
which is equivalent to 
\begin{equation}\label{nonnn}
\begin{aligned}
& \int |N^3v(N(x-y))| \big\|\dxx\lppt(t,x,y)\big\|_{L^{2}(dt)L^{2}(d(x+y))} d(x-y)\\
 &\le  C\e\|\lppt\|_{\mathcal{N}_1(\Lambda)}.
\end{aligned}
\end{equation}
To see this, note that if for $\lppt$, we have $|\xi-\eta|<20N$, $|\xi|<20N$ or $|\eta|<20N$,  then the left side of \eqref{nonnn} is easily controlled by 
\begin{equation}
\begin{aligned}
& \int |N^3v(N(x-y))| \big\|\dxx\lppt(t,x,y)\big\|_{L^{2}(dt)L^{2}(d(x+y))} d(x-y)\\
 &\le  C\e\|\lppt\|_{\lc}.
\end{aligned}
\end{equation}
Thus if we denote $P_{>N}=P_{|\xi-\eta|\ge20N}P_{|\xi|\ge 20N}P_{|\eta|\ge 20N}$, we are further reduced to showing that 
\begin{equation}\label{nonnn1}
\begin{aligned}
& \int |N^3v(N(x-y))| \big\|\dxx P_{>N}\lppt(t,x,y)\big\|_{L^{2}(dt)L^{2}(d(x+y))} d(x-y)\\
 &\le  C\e\|\lppt\|_{\mathcal{N}_1(\Lambda)}.
\end{aligned}
\end{equation}
As in the proof of main theorem for the linear equation, we shall divide our discussion into two cases.

Case 1: $|\xi+\eta|< \frac{N}{10}$.

In this case, by H\"older's inequality,
\begin{align*}
& \int |N^3v(N(x-y))| \big\|\dxx P_{>N}\lppt(t,x,y)\big\|_{L^{2}(dt)L^{2}(d(x+y))} d(x-y)\\
 &\le  C\e N^{\frac12} \big\|\dxx P_{>N}P_{|\xi+\eta|< \frac{N}{10}}\lppt(t,x,y)\big\|_{L^6(d(x-y))L^{2}(dt)L^{2}(d(x+y))}\\
  &\le  C\e N^{\frac12} \big\|\dxx P_{>N}P_{|\xi+\eta|< \frac{N}{10}}\lppt(t,x,y)\big\|_{L^{2}(dt)L^6(d(x-y))L^{2}(d(x+y))}\\
    &\le  C\e  \big\|\dxx\dyy P_{|\xi-\eta|\ge10N}P_{|\xi+\eta|< \frac{N}{10}}\lppt(t,x,y)\big\|_{L^{2}(dt)L^6(d(x-y))L^{2}(d(x+y))}\\
     &\le  C\e\|\Lm_p^T\|_{\mathcal{N}_1(\Lambda)},
\end{align*}
where in the third inequality we used Bernstein's  inequality  along with the fact that $|\xi-\eta|\ge 20N$.

Case 2: $|\xi+\eta|\ge \frac{N}{10}$.

In this case, by H\"older's inequality,
\begin{align*}
& \int |N^3v(N(x-y))| \big\|\dxx P_{>N}\lppt(t,x,y)\big\|_{L^{2}(dt)L^{2}(d(x+y))} d(x-y)\\
 &\le  C\e N^{\frac12} \big\|\dxx P_{>N}P_{|\xi+\eta|\ge \frac{N}{10}}\lppt(t,x,y)\big\|_{L^6(d(x-y))L^{2}(dt)L^{2}(d(x+y))}\\
  &\le  C\e N^{\frac12} \big\|\dx P_{>N}P_{|\xi+\eta|\ge \frac{N}{10}}\lppt(t,x,y)\big\|_{L^{2}(dt)L^6(d(x-y))L^{2}(d(x+y))}\\
    &\quad+  C\e N^{\frac12} \big\|\dy P_{>N}P_{|\xi+\eta|\ge \frac{N}{10}}\lppt(t,x,y)\big\|_{L^{2}(dt)L^6(d(x-y))L^{2}(d(x+y))}\\
    &\le  C\e  \big\|\dx\dy P_{|\xi|\ge10N}P_{|\eta|\ge {10N}}P_{|\xi+\eta|\ge \frac{N}{10}}\lppt(t,x,y)\big\|_{L^{2}(dt)L^6(d(x-y))L^{2}(d(x+y))}\\
     &\le  C\e\|\lppt\|_{\mathcal{N}_1(\Lambda)}
\end{align*}
where in the third inequality we used Bernstein's  inequality in rotated coordinates along with the fact that $|\xi|\ge 20N$ and $|\eta|\ge 20N$.

\end{proof}
We continue with estimates for $\|\lcct\|_{\mathcal{N}_2(\Lambda)}$. This is an easy version of the previous theorem. Using Lemma~\ref{xco1}-Lemma~\ref{xco3} and Strichartz 
estimates (Theorem~\ref{2.1'}), and then applying Lemma~\ref{nonlinear}-\ref{nonlinear2} in this section to handle the nonlinear terms, we get
\begin{theorem}\label{lcb}
Let $[0,T]$ be as above, there exist a universal constant $C$ such that
\begin{equation}
\begin{aligned}
\|\lcct\|&_{\mathcal{N}_2(\Lambda)}\le C\|\dx\dy\Lambda_c(0, \cdot)\|_{L^2}+ C\e \|\lcct\|_{\mathcal{N}_2(\Lambda)} \\
&\quad\qquad +C\e \|\lppt\|_{\mathcal{N}_1(\Lambda)}\|\gct\|_{\mathcal{S}_{x,y}}.
\end{aligned} 
\end{equation}
\end{theorem}
Using Strichartz estimates for $\S_\pm$ and Lemma~\ref{nonlinear}-\ref{nonlinear2}, we get
\begin{theorem}\label{gcb}
Let $[0,T]$ be as above, there exist a universal constant $C$ such that
 \begin{equation}
\begin{aligned}
\|\gct\|&_{\mathcal{S}_{x,y}}\le C\|\dx\dy\gc(0, \cdot)\|_{L^2}+ C\e\|\gct\|_{\mathcal{S}_{x,y}}\\
&\quad\qquad +C\e \|\lppt\|_{\mathcal{N}_1(\Lambda)}\|\lcct\|_{\mathcal{N}_2(\Lambda)}
\end{aligned}\ed 
 \begin{equation}
\begin{aligned}
\|\gpt\|&_{\mathcal{S}_{x,y}}\le C\|\dx\dy\gp(0, \cdot)\|_{L^2}+ C\e\|\gpt\|_{\mathcal{S}_{x,y}}\\
&\quad\qquad +C\e \|\lppt\|_{\mathcal{N}_1(\Lambda)}+C\e \|\lppt\|_{\mathcal{N}_1(\Lambda)}\|\lppt\|_{\mathcal{N}_1(\Lambda)}.
\end{aligned} \ed
\end{theorem}
For later use, let us denote 
 \begin{equation}
\begin{aligned}
X(T)=& \|\lcct\|_{\mathcal{N}_2(\Lambda)}+\|\gct\|_{\mathcal{S}_{x,y}}\\
Y(T)=& \|\lppt\|_{\mathcal{N}_1(\Lambda)}+\|\gpt\|_{\mathcal{S}_{x,y}}.
\end{aligned} \ed
We want to show that $X(T), \,Y(T)$ depends continuously on $T$. To see this, for any fixed $T\ge0$,  we have 
$$\Big|\|\lcct\|_{\mathcal{N}_2(\Lambda)}-\|\Lm_{c,T+\delta}\|_{\mathcal{N}_2(\Lambda)}\Big|\ls \|\Lm_{c,T+\delta}-\lcct\|_{\mathcal{N}_2(\Lambda)}.
$$
And note that $\Lm_{c,T+\delta}-\lcct$ satisfy
$$\S(\Lm_{c,T+\delta}-\lcct)=\chi_{[T, T+\delta]}\Big(-\{V_N*\rho, \lcc\}-\big((V_N\bar\Gamma_p)\circ\lcc+(V_N\lpp)\circ\gc\big)_{symm} \Big),$$
with 0 initial condition. It is not hard to see that, by crude energy estimates
 \begin{equation}
 \|\Lm_{c,T+\delta}-\lcct\|_{\mathcal{N}_2(\Lambda)}\le \delta\, C(T,N) \ed
for some constant $C(T,N)$ that depend on $T, N$, which implies the continuity of $\|\lcct\|_{\mathcal{N}_2(\Lambda)}$. The continuity of other norms in $X(T)$ and 
$Y(T)$ can be proved similarly.

We can now state and prove the main theorem of this section
\begin{theorem}\label{nonlinear main}
Assume $\Lm$, $\Gamma$ and $\phi$ are smooth solutions to the HFB system, with finite energy per particle, uniformly in $N$ (see \eqref{initial condition}), which implies 
\eqref{lm 1233}-\eqref{lm 123''}, we have
 \begin{equation}
\begin{aligned}
& \|\Lm_c\|_{\mathcal{N}_2(\Lambda)}+\|\gc\|_{\mathcal{S}_{x,y}} \le C \\
& \|\Lm_p\|_{\mathcal{N}_1(\Lambda)}+\|\gp\|_{\mathcal{S}_{x,y}} \le C.
\end{aligned}\ed
\end{theorem}
\begin{proof}
By Theorem~\ref{lpb}, \ref{lcb}, \ref{gcb}, and the size of initial conditions \eqref{lm 1233}-\eqref{lm 123''}, we have 
\begin{equation}\label{maxs}
\begin{aligned}
X(T)\le & CC_0 +C\e X(T)Y(T) \\
Y(T)\le &CC_0+C \e Y(T)^2+C\e X(T).
\end{aligned}\ed
At this stage, we will need to assume $\e$ is small, where the smallness may depend on $C_0$. 
Without loss of generality, let's assume $CC_0=1$ and $C\e\le \frac{1}{10}$, then we have the following simplified version of \eqref{maxs}
\begin{equation}\label{maxs2}
\begin{aligned}
X(T)\le & 1 +\frac{1}{10} X(T)Y(T) \\
Y(T)\le &1+ \frac{1}{10} Y(T)^2+\frac{1}{10} X(T).
\end{aligned}\ed
If $Y(T)\le 4$, by the first line, we have $X(T)\le 2$, and if we plug this into the second line, we get $Y(T)\le 3$.
By continuity, since $Y(0)\le 1$, we always have $Y(T)\le 3$, and thus $X(T)\le 2$, for all $T\ge0$, which concludes our proof.
\end{proof}

\newsection{Estimates for $\mbox{\rm sh} (2k)$, $p_2=\overline{\mbox{\rm sh} (k)}\circ \sh$ and $\sh$}

The equations for $\sht=N\lpp$ and $p_2=N\gp$ are
\begin{align*}
&\S\, \sht+ \{V_N*\rho, \sht\}+\big((V_N\bar\Gamma)\circ\sht+(V_N\Lm)\circ p_2\big)_{symm}=-\frac{V_N}{2}\Lm\\
&\S_{\pm}\,p_2+[V_N*\rho, p_2]+\big((V_N\Gamma)\circ p_2+(V_N\bar\Lm)\circ \sht\big)_{skew}=0.
\end{align*}

To handle the inhomogeneous term $-\frac{V_N}{2}\Lm$, we shall need the following lemma.

\begin{lemma}\label{n3}
Let $\S u=-\frac{V_N}{2}\Lm_p$ with $u(0,\cdot)=0$, we have 
$$\|u\|_{\mathcal{S}_{x,y}}\ls \|\Lm_p\|_{\mathcal{N}_1}
$$
where $\|\cdot\|_{\mathcal{N}_1}$ is defined as in \eqref{lambdap norm}.
\end{lemma}
Note that the above result also hold if one replace $\Lm_p$ by $\Lm_c$ and replace $\|\Lm_p\|_{\mathcal{N}_1}$ by $\|\Lm_c\|_{\mathcal{N}_2}$, which is a direct consequence of 
Theorem~\ref{Spsi} and H\"older's inequality.
\begin{proof}
First note that if for $\Lm$, we have $|\xi-\eta|<20N$, $|\xi|<20N$ or $|\eta|<20N$,  then by Theorem~\ref{Spsi} and H\"older's inequality, we have
\begin{equation}\nonumber
\begin{aligned}
\|u\|_{\mathcal{S}_{x,y}}\le  C\e\|\Lm_p\|_{\lc}.
\end{aligned}
\end{equation}
Thus if we denote $P_{>N}=P_{|\xi-\eta|\ge20N}P_{|\xi|\ge 20N}P_{|\eta|\ge 20N}$, if suffices to prove Lemma~\ref{n3} with $\Lm$ replaced by $P_{>N}\Lm$.
As in the proof of main theorem for the linear equation, we shall divide our discussion into two cases.

Case 1: $|\xi+\eta|< \frac{N}{10}$.

In this case, by Theorem~\ref{Spsi} and H\"older's inequality, we have 
 \begin{equation}
\begin{aligned}
&\|u\|_{\mathcal{S}_{x,y}}\\
 &\le C  \big\|V_N\dxx P_{>N}P_{|\xi+\eta|< \frac{N}{10}}\Lm_p\big\|_{L^1(d(x-y))L^{2}(dt)L^{2}(d(x+y))}\\
&\quad+C  \big\|V_N\dt P_{>N}P_{|\xi+\eta|< \frac{N}{10}}\Lm_p\big\|_{L^1(d(x-y))L^{2}(dt)L^{2}(d(x+y))}\\
&\le C\e N^{\frac12} \big\|\dxx P_{>N}P_{|\xi+\eta|< \frac{N}{10}}\Lm_p\big\|_{L^6(d(x-y))L^{2}(dt)L^{2}(d(x+y))}\\
&\quad+C\e N^{\frac12} \big\|\dt P_{>N}P_{|\xi+\eta|< \frac{N}{10}}\Lm_p\big\|_{L^6(d(x-y))L^{2}(dt)L^{2}(d(x+y))}\\
&\le C\e N^{\frac12} \big\|\dxx P_{>N}P_{|\xi+\eta|< \frac{N}{10}}\Lm_p\big\|_{L^{2}(dt)L^6(d(x-y))L^{2}(d(x+y))}\\
&\quad+C\e N^{\frac12} \big\|\dt P_{>N}P_{|\xi+\eta|< \frac{N}{10}}\Lm_p\big\|_{L^{2}(dt)L^6(d(x-y))L^{2}(d(x+y))}\\
  &\le  C\e\big\|\dt\dyy P_{|\xi-\eta|\ge10N}P_{|\xi+\eta|< \frac{N}{10}}\Lm_p\big\|_{L^{2}(dt)L^6(d(x-y))L^{2}(d(x+y))}\\
    &\quad+  C\e  \big\|\dxx\dyy P_{|\xi-\eta|\ge10N}P_{|\xi+\eta|< \frac{N}{10}}\Lm_p\big\|_{L^{2}(dt)L^6(d(x-y))L^{2}(d(x+y))}\\
     &\le  C\e\|\Lm_p\|_{\mathcal{N}_1(\Lambda)}
\end{aligned}\ed
where in the fourth inequality we used Bernstein's  inequality  along with the fact that $|\xi-\eta|\ge 20N$.

Case 2: $|\xi+\eta|\ge \frac{N}{10}$.

In this case, by Strichartz (Theorem~\ref{2.1'}) and H\"older's inequality,
 \begin{equation}
\begin{aligned}
&\|u\|_{\mathcal{S}_{x,y}}\\
 &\le C  \big\|V_N P_{>N}P_{|\xi+\eta|>\frac{N}{10}}\Lm_p\big\|_{L^{2}(dt)L^{\frac65}(d(x-y))L^{2}(d(x+y))}\\
&\le C\e N\big\| P_{>N}P_{|\xi+\eta|\ge \frac{N}{10}}\Lm_p\big\|_{L^{2}(dt)L^6(d(x-y))L^{2}(d(x+y))}\\
  &\le  C\e\big\|\dx\dy P_{|\xi|\ge10N}P_{|\eta|\ge10N}P_{|\xi+\eta|< \frac{N}{10}}\Lm_p\big\|_{L^{2}(dt)L^6(d(x-y))L^{2}(d(x+y))}\\
     &\le  C\e\|\Lm_p\|_{\mathcal{N}_1(\Lambda)}
\end{aligned}\ed
where in the third inequality we used Bernstein's  inequality in rotated coordinates along with the fact that $|\xi|\ge 20N$ and $|\eta|\ge 20N$.
\end{proof}

Now we shall estimate the other nonlinear terms in dual Strichartz norms
\begin{equation}\label{shk1}
\begin{aligned} 
\|&\big(V_N*\rho(t,x)\big) \sht(t,x,y)\|_{L^{\frac85}(dt)L^{\frac43}(dx)L^2(dy)}\\
&\qquad+\|\big((V_N\bar\Gamma)\circ\sht\|_{L^{\frac85}(dt)L^{\frac43}(dx)L^2(dy)}\\
&\le C\e \sup_z \|\Gamma (t, x+z, x)\|_{L^{8}(dt)L^{\frac{12}{7}}(dx)}\|\sht\|_{L^{2}(dt)L^{6}(dx)L^2(dy)}\\
&\le C\e \|\sht\|_{L^{2}(dt)L^{6}(dx)L^2(dy)},
\end{aligned}
\end{equation}
where we used Sobolev and Lemma~\ref{a priori} in the last inequality. The above estimate still hold if we replace $\bar\Gamma$ by $\bar\Lm_c$. Similarly 
\begin{equation}\label{shk2}
\begin{aligned}
\|&\big(V_N*\rho(t,x)\big) p_2(t,x,y)\|_{L^{\frac85}(dt)L^{\frac43}(dx)L^2(dy)}\\
&\qquad+\|\big((V_N\Gamma)\circ p_2\|_{L^{\frac85}(dt)L^{\frac43}(dx)L^2(dy)}\\
&\le C\e \|p_2\|_{L^{2}(dt)L^{6}(dx)L^2(dy)}.
\end{aligned}
\end{equation}
And the above estimate still hold if we replace $\Gamma$ by $\Lm_c$. 
Also, since for fixed $z$, the following holds, uniformly in $z$
\begin{equation}\label{shk3}
\begin{aligned}
&\|\Lm(t,x,x+z)p_2(t,x+z,y)\|_{L^{\frac43}(dt)L^{\frac32}(dx)L^2(dy)} \\
&\le \|\Lm(t,x,x+z)\|_{L^{2}(dt)L^{3}(dx)} \|p_2(t,x,y)\|_{L^{4}(dt)L^{3}(dx)L^2(dy)} \\
&\le C\|\dx\Lm(t,x,x+z)\|_{L^{2}(dt)L^{2}(dx)} \|p_2(t,x,y)\|_{L^{4}(dt)L^{3}(dx)L^2(dy)}.
\end{aligned}
\end{equation}
If we repeat the argument in the proof of Lemma~\ref{nonlinear2} and using \eqref{nonnn}, we have 
\begin{align*}
\|&\big(V_N\Lm_p\big)\circ p_2\|_{L^{\frac85}(dt)L^{\frac43}(dx)L^2(dy)}\\
&\le C \int |N^3v(Nz)| \|\dx\lpp(t,x,x+z)\|_{L^{2}(dt)L^{2}(dx)}\|p_2\|_{L^{4}(dt)L^{3}(dx)L^2(dy)} dz \\
 &\le  C\e\|\Lm\|_{\mathcal{N}_1(\Lambda)}\|p_2\|_{L^{4}(dt)L^{3}(dx)L^2(dy)}.
\end{align*}
Similarly, we also have 
\begin{align*}
\|&\big(V_N\bar\Lm_p\big)\circ \sht\|_{L^{\frac85}(dt)L^{\frac43}(dx)L^2(dy)}\\
&\le C \int |N^3v(Nz)| \|\dx\lpp(t,x,x+z)\|_{L^{2}(dt)L^{2}(dx)}\|\sht\|_{L^{4}(dt)L^{3}(dx)L^2(dy)} dz \\
 &\le  C\e\|\Lm\|_{\mathcal{N}_1(\Lambda)}\|\sht\|_{L^{4}(dt)L^{3}(dx)L^2(dy)}.
\end{align*}
If we choose $\e$ small enough such that $\e C\le \frac{1}{10}$, combing the above estimates, and using the fact that $\|\Lm_p\|_{\mathcal{N}_1(\Lambda)}+\|\Lm_c\|_{\mathcal{N}_2(\Lambda)}\ls 1$ from the main theorem of the last section, we get, by Strichartz 
 \begin{equation}
\begin{aligned}
\|& \sht\|_{\mathcal{S}_{x,y}}+\|p_2\|_{\mathcal{S}_{x,y}}\\
&\le C\Big(\| \sht(0,\cdot)\|_{L^2}+\|p_2(0,\cdot)\|_{L^2}\Big) +C.\end{aligned}\ed

\newsection{Estimates for the condensate $\phi$}

The non-linear equation for $\phi$ can be regarded as a linear equation on a background given by $\Gamma$ and $\Lm$, for which we already have estimates:
\bgd\ad
\Big\{&\frac1i\partial_t-\Delta_{x_1}\Big\}\phi(x_1)\\
&=-\int dy \{V_N(x_1-y)\Gamma(y,y)\}\phi(x_1)  \\
&\,\,\,\,-\int dy \{V_N(x_1-y)\Gamma_p(y,x_1)\}\phi(y)  \\
&\,\,\,\,+\int dy \{V_N(x_1-y)\lpp(x_1,y)\}\bar\phi(y) .
\ead\ed
Define the standard Strichartz spaces
\begin{align}
\|\phi\|_{\mathcal{S}}
= \sup_{  (p,q) \, \, admissible}\|\phi\|_{L^{p}(dt) L^{q}(dx)}.\notag
\end{align}
\begin{proof}
We shall estimate the right hand side of the equation for $\phi$ in dual Strichartz norms, if we repeat the proof of Lemma~\ref{nonlinear}-\ref{nonlinear2}, 
it is not hard to show 
\bgd\ad
\|\langle \nabla\rangle^{\frac12}\int dy \{V_N(x_1-y)\Gamma(y,y)\}\phi(x_1) \|_{L^{\frac85}(dt)L^{\frac43}(dx)} \le C\e \|\langle \nabla\rangle^{\frac12}\phi\|_{L^{2}(dt)L^{6}(dx)}  \\
\|\langle \nabla\rangle^{\frac12}\int dy \{V_N(x_1-y)\Gamma_p(y,x_1)\}\phi(y) \|_{L^{\frac85}(dt)L^{\frac43}(dx)} \le C\e \|\langle \nabla\rangle^{\frac12}\phi\|_{L^{2}(dt)L^{6}(dx)}  \\
\|\langle \nabla\rangle^{\frac12}\int dy \{V_N(x_1-y)\lpp(x_1,y)\}\bar\phi(y)\|_{L^{\frac43}(dt)L^{\frac32}(dx)} \le C\e \|\langle \nabla\rangle^{\frac12}\phi\|_{L^{4}(dt)L^{3}(dx)}.
\ead\ed
Thus, 
\begin{align*}
\| \langle \nabla\rangle^{\frac12}\phi\|_{\mathcal{S}}
\le C\|\langle \nabla\rangle^{\frac12} \phi (0,\cdot)\|_{L^2} +3C\e \| \langle \nabla\rangle^{\frac12}\phi\|_{\mathcal{S}}.
\end{align*}
This gives us desired result by taking $\e C\le \frac{1}{10}$.
\end{proof}
\newsection{Remarks on the collapsing norm of $\Lm$.}

In this section, we shall see how we can use Theorem~\ref{main} to prove the following
\begin{theorem}
Let $\Lambda$ satisfy
\begin{equation}\label{f.1}
\begin{aligned}
\S\Lambda(t,x,y)&+N^2v\big(N(x-y)\big)\Lambda(t,x,y)=G(t,x,y), \,\, \Lambda(0,\cdot)=\Lambda_0
\end{aligned} 
\end{equation}
where we assume $v$ satisfy \eqref{vcondition}, we have
\begin{equation}\label{f.2}
\begin{aligned}
 \|\dxx \Lm&\|_{collapsing} + \|\dt \Lm\|_{collapsing} \\
& \ls \|\dx\dy G\|_{\mathcal{S}'_{r}} +\|\dx\dy \Lm_0\|_{L^2}.
\end{aligned}
\end{equation}
\end{theorem}
Recall that we already have desired estimates if we replace $\|\cdot\|_{collapsing}$ by  $\|\cdot\|_{\lc}$,  by using the results in Theorem~\ref{main}.

\begin{proof}
We shall first treat the homogeneous equation, let 
\begin{equation}\label{f.3}
\begin{aligned}
\S\Lambda(t,x,y)&+N^2v\big(N(x-y)\big)\Lambda(t,x,y)=0 \\
 \Lambda(0,\cdot)&=\Lambda_0.
\end{aligned} 
\end{equation}
Let $H=-\Delta_x-\Delta_y + N^2v(N(x-y))$, $H_0=-\Delta_x-\Delta_y$. If we let $W$ denote the wave operator in Yajima's paper \cite{yajima},  acting in $x-y$ direction, we have
$e^{i t H}= We^{i t H_0}W^*$, where $^*$ denotes the dual operator. We also have, for the potential $v$ satisfying \eqref{vcondition}, $W$ is a bounded operator from $L^p\rightarrow L^p$ for any $1\le p\le \infty$, with bound independent of $N$ (see \cite{erdHos2009rigorous} Proposition 5.1). Moreover, by calculating the integral kernel of $W$ explicitly,  the $L^p\rightarrow L^p$
boundness of $W$ extends to the space of $L^2$ valued function, ).  
Thus, 
\begin{align*}
&\sup_{x-y}\|\dxx e^{-i t H}\Lambda_0 \|_{L^2(dt d(x+y) )}
+ \sup_{x-y}\|||\frac{\partial}{\partial t}|^{1/4} e^{-i t H}\Lambda_0 \|_{L^2(dt d(x+y) )}\\
&=\sup_{x-y}\|W\dxx e^{i t H_0}W^*\Lambda_0 \|_{L^2(dt d(x+y) )}
+ \sup_{x-y}\|W|\frac{\partial}{\partial t}|^{1/4} e^{i t H_0}W^*\Lambda_0 \|_{L^2(dt d(x+y) )}\\
&\lesssim \sup_{x-y}\|\dxx e^{i t H_0}W^*\Lambda_0 \|_{L^2(dt d(x+y) )}
+ \sup_{x-y}\||\frac{\partial}{\partial t}|^{1/4} e^{i t H_0}W^*\Lambda_0 \|_{L^2(dt d(x+y) )}\\
&\lesssim \|\dx\dy W^*\Lambda_0\|_{L^2(dxdy)},
\end{align*}
where we used Lemma~\ref{xco2}-\ref{xco3} in the last inequality.

Recall $W^*=\lim_{t \to \infty} e^{i t H_0}e^{-i t H}$, where the limit exists in the strong operator topology, see e.g., Proposition 5.1 in \cite{erdHos2009rigorous} and also \cite{georgiev2005existence} for a proof of the existence of strong limit when  $v$ satisfy \eqref{vcondition}. Thus, for each fixed $g\in C_0^\infty(\R^6)$ with $\|g\|_2=1$
\bgd
\ad
& \langle \dx\dy W^*\Lambda_0, g\rangle\\
&=\langle W^*\Lambda_0, \dx\dy g\rangle\\
&= \lim_{t \to \infty} \langle e^{i t H_0}e^{-i t H} \Lambda_0, \dx\dy g\rangle\\
&= \lim_{t \to \infty} \langle e^{-i t H} \Lambda_0, \dx\dy e^{-i t H_0} g\rangle\\
&= \lim_{t \to \infty} \langle \dx\dy e^{-i t H} \Lambda_0,  e^{-i t H_0} g\rangle\\
&\le  \|\dx\dy e^{-i t H} \Lambda_0\|_{L^\infty (dt)L^2(dxdy)} \\ 
&\le  \|\dx\dy e^{-i t H} \Lambda_0\|_{\mathcal{S}_{x,y}} \\ 
&\le \|\dx\dy \Lm_0\|_{L^2}
\ead
\ed
where the second equality is a consequence of the existence of strong limit, and  in the last inequality we used the special case of Theorem~\ref{main} with $G=H=0$. By taking supremum among all choices of $g$, we have 
$$\|\dx\dy W^*\Lambda_0\|_{L^2(dxdy)}\ls \|\dx\dy \Lm_0\|_{L^2}.
$$

Now we shall consider the inhomogeneous equation, let
\begin{equation}\label{f.4}
\begin{aligned}
\S\Lambda(t,x,y)&+N^2v\big(N(x-y)\big)\Lambda(t,x,y)=G(t,x,y),\,\, \Lambda(0,\cdot)=0.
\end{aligned} 
\end{equation}
If we let $H$ be defined as above, we have by Duhamel's formula
$$\Lm=\int_0^t e^{-i(t-s)H}G(s,\cdot)ds.$$
We shall first show that 
\begin{align}
\sup_{x-y}\|&\dxx\int_0^t e^{-i(t-s)H}G(s,\cdot)ds \|_{L^2(dt d(x+y) )} \nonumber
\\ &\lesssim \|\dx\dy G\|_{\sr}.\label{fug}
\end{align}
To prove \eqref{fug}, it suffices to show that, for any fixed $x-y$ and fixed $T>0$,
\begin{align*}
\|&\chi_{[0,T]}(t)\dxx\int_0^t e^{-i(t-s)H}G(s,\cdot)ds \|_{L^2(dt d(x+y) )} \nonumber
\\ &\lesssim \|\dx\dy G\|_{\sr}
\end{align*}
which, by the Christ--Kiselev lemma, is a consequence of 
\begin{equation} \label{fug1}
\begin{aligned}
\|&\dxx\int_0^T e^{-i(t-s)H}G(s,\cdot)ds \|_{L^2(dt d(x+y) )} 
\\ &\lesssim \|\dx\dy G\|_{\sr}.
\end{aligned}
\end{equation}
If we apply the homogeneous estimates proved above, the left side of \eqref{fug1} is bounded by
\begin{equation} \label{fug2}
\begin{aligned}
\|&\dx\dy\int_0^T e^{isH}G(s,\cdot)ds \|_{L^2(dxdy)}  \\
&=\|\dx\dy e^{iTH}\int_0^T e^{-i(T-s)H}G(s,\cdot)ds \|_{L^2(dxdy)} \\
&\ls \|\dx\dy\int_0^T  e^{-i(T-s)H}G(s,\cdot)ds \|_{L^2(dxdy)} \\
&\lesssim \|\dx\dy G\|_{\sr},
\end{aligned}
\end{equation}
where in the third line we used the special case of Theorem~\ref{main} with $G=H=0$ (estimates for the homogeneous equation), and in the 
last line we used the  special case of Theorem~\ref{main} with $\Lm_0=H=0$(inhomogeneous equation with zero initial data).

To prove inhomogeneous estimate for $\dt$ derivative, as before we can not use the Christ--Kiselev lemma, we shall follow the ideas in the prove of Lemma~\ref{xco3}.
Write $\Lambda=\Lambda_1+\Lambda_2$, where
\begin{align*}
 \S \, \Lambda_1&= G(t,x,y), \,\,\,\,\text{with initial conditions 0}\\
 \S \, \Lambda_2&= -N^2v\big(N(x-y)\big)\Lambda(t,x,y), \,\,\,\,\text{with initial conditions 0}.
\end{align*}
For both $\Lambda_1$ and $\Lambda_2$, it suffices to consider the region where 
\begin{equation}\label{recondi}
 \tau^{\frac12} \ge 10(1+|\xi+\eta|),  |\xi|\ge 20N, \,\,  |\eta|\ge 20N\,\,\text{and}\, \, |\xi-\eta|\ge 20N,
\end{equation}
since otherwise we have $\|\dt\Lm\|_{collapsing}\ls \|\dxx\Lm\|_{collapsing}$, or we already know how to control collapsing norm when $|\xi|< 20N$, $|\eta|< 20N$, or  $|\xi-\eta|< 20N$.

The estimate for $\Lm_1$ just follows from Lemma~\ref{xco3} directly, and for $\Lm_2$, we claim that it suffices to show 
\begin{proposition}\label{8.2}
Let $\S u=f$ with $u(0,\cdot)=0$, then if the Fourier support $\tau, \xi, \eta$ of $u$ satisfy \eqref{recondi}, we have 
\bgd\ad\|&\dt u\|_{collapsing} 
 \lesssim  \min\Big\{\|\dx\dy f\|_{L^{2}(dt)L^{\frac65}(d(x-y))L^{2}(d(x+y))}, \\ 
&\|\dt\dyy f\|_{L^{2}(dt)L^{\frac65}(d(x-y))L^{2}(d(x+y))}+N^{-\frac12}\|\dt f\|_{L^{2}(dtdxdy)} \Big\}.
\ead\ed\end{proposition}
\begin{remark}
Note that the $\tau$ support of $u$ may be different from $f$, so $\tau$ support of $f$ may not satisfy \eqref{recondi} , but the  $\xi, \eta$ support of $f$ does satisfy \eqref{recondi}. 
\end{remark}
We shall first see how we can apply the Proposition to get desired results.

Case 1: $|\xi+\eta|\ge\frac{N}{10}$.

In this case, 
$$ \S \,P_{|\xi|\ge 20N}P_{|\eta|\ge 20N} \Lambda_2\sim -N^2v\big(N(x-y)\big)P_{|\xi|\ge 19N}P_{|\eta|\ge 19N} \Lambda(t,x,y). $$
Thus, by Proposition~\ref{8.2}
\begin{align*}
\|&\dt \Lm_2\|_{collapsing} \\
 \lesssim & \|\dx\dy N^2v\big(N(x-y)\big)P_{|\xi|\ge 19N}P_{|\eta|\ge 19N} \Lambda(t,x,y)\|_{L^{2}(dt)L^{\frac65}(d(x-y))L^{2}(d(x+y))}\\
  \lesssim & \|\dx\dy P_{|\xi|\ge 19N}P_{|\eta|\ge 19N} \Lambda(t,x,y)\|_{L^{2}(dt)L^{6}(d(x-y))L^{2}(d(x+y))},
 \end{align*}
where in the last inequality we used the fact that when $|\xi|\ge 19N$, $|\eta|\ge 19N$, $\dx\dy$ essentially only fall on $\Lm$, as well as H\"older's inequality.
Recall that when $|\xi+\eta|>\frac{N}{10}$, the last line appears as part of the norm for $\Lm$ in the proof of Theorem~\ref{main} (see \eqref{1655}), thus it is bounded by
$\|\dx\dy G\|_{\sr}
$ for $\Lm$ satisfying \eqref{f.4}.

Case 2: $|\xi+\eta|< \frac{N}{10}$.

In this case, 
$$ \S \,P_{|\xi+\eta|< \frac{N}{10}}P_{|\xi-\eta|\ge 20N} \Lambda_2\sim -N^2v\big(N(x-y)\big)P_{|\xi+\eta|< \frac{N}{10}}P_{|\xi-\eta|\ge 19N} \Lambda(t,x,y). $$
Thus, by Proposition~\ref{8.2}
\begin{align*}
\|&\dt \Lm_2\|_{collapsing} \\
 \lesssim & \|\dt\dyy N^2v\big(N(x-y)\big)P_{|\xi+\eta|< \frac{N}{10}}P_{|\xi-\eta|\ge 19N}\Lambda(t,x,y)\|_{L^{2}(dt)L^{\frac65}(d(x-y))L^{2}(d(x+y))}\\
  &+ \|\dt N^{\frac32}v\big(N(x-y)\big)P_{|\xi+\eta|< \frac{N}{10}}P_{|\xi-\eta|\ge 19N}\Lambda(t,x,y)\|_{L^{2}(dtdxdy)}\\
  \lesssim & \|\dt\dyy P_{|\xi+\eta|< \frac{N}{10}}P_{|\xi-\eta|\ge 19N} \Lambda(t,x,y)\|_{L^{2}(dt)L^{6}(d(x-y))L^{2}(d(x+y))} \\
  &+ N^{\frac12}\|\dt  P_{|\xi+\eta|< \frac{N}{10}}P_{|\xi-\eta|\ge 19N} \Lambda(t,x,y)\|_{L^{2}(dt)L^{6}(d(x-y))L^{2}(d(x+y))}  \\
   \lesssim & \|\dt\dyy P_{|\xi+\eta|< \frac{N}{10}}P_{|\xi-\eta|\ge 19N} \Lambda(t,x,y)\|_{L^{2}(dt)L^{6}(d(x-y))L^{2}(d(x+y))}
 \end{align*}
where in the second inequality we used the fact that when $|\xi-\eta|\ge 19N$,  $\dyy$ essentially only fall on $\Lm$, as well as H\"older's inequality. And 
in the last inequality we used Bernstein's inequality and the fact that $|\xi-\eta|\ge 20N$.

Recall that when $|\xi+\eta|<\frac{N}{10}$, the last line appears as part of the norm for $\Lm$ in the proof of Theorem~\ref{main} (see \eqref{867} ), thus it is bounded by
$\|\dx\dy G\|_{\sr}
$ for $\Lm$ satisfying \eqref{f.4}.

Thus it remains to prove Proposition~\ref{8.2}, to see this, we shall follow the ideas in the proof of Lemma~\ref{xco3}.

\noindent{$ Proof \,\,of \,\,Proposition\,\,$~\ref{8.2}}:

Case 1: If $|\tau|^{\frac{1}{2}} > 2 (|\xi|+ |\eta|)$.

Write $u=u^1+u^2$, where
\begin{equation}\label{8123'}
\begin{aligned}
\mathcal{F} \, u^1&=\frac{\mathcal{F}f}{\tau+|\xi|^2+|\eta|^2}, \,\,\text{this no longer has initial conditions 0} \\
 \S \, u^2&= 0, \,\,\,\,\text{a correction so that $u^1+u^2$ has initial condition 0}.
 \end{aligned}
 \end{equation}
In this case, it suffices to control $u_1$ since $u_2$ is only supported where $|\tau|=|\xi|^2+|\eta|^2$. The goodness about $u_1$ is that it has the same Fourier support with $f$.
The strategy is based on
\begin{align}
\|\dt u_1 \|&_{L^{\infty}(d(x-y)) L^2(d(x+y)dt)}
=\|\tau^{\frac14}\frac{\mathcal{F}f}{\tau+|\xi|^2+|\eta|^2} \|_{L^{\infty}(d(x-y))L^2(d \tau d (\xi+\eta))} \nonumber\\
&\lesssim \|\int|\tau|^{\frac14} \Big|\frac{\mathcal{F}\Big(f\Big)}{\tau+|\xi|^2+|\eta|^2}\Big|d(\xi-\eta)\|_{L^2(d \tau d(\xi+\eta))}. \label{8'}
\end{align}

By Cauchy-Schwarz, we have
\bgd\ad
&\mathrm{RHS} \eqref{8'}
\lesssim \|\frac{|\tau|^{\frac{1}{4}}}{|\tau|}\int_{|\xi-\eta|< |\tau|^{\frac{1}{2}}} |\mathcal{F}f| d(\xi-\eta)\|_{L^2(d \tau d(\xi+\eta))}\\
&\lesssim A \big\|{|\nabla_y|^{\frac12}}{|\nabla_x|}^{-\frac 12}f\big\|_{L^2(dtd(x-y)d(x+y))}
\ead\ed where
\bgd\ad A= \sup_{\tau, \xi+\eta}\frac{|\tau|^{\frac{1}{4}}}{|\tau|} \left(\int_{|\xi-\eta|< |\tau|^{\frac{1}{2}}}
\frac{|\xi|}{|\eta|}d(\xi-\eta)\right)^{\frac{1}{2}}.
\ead\ed Changing variables, this is something like
\begin{align*}
A= \sup_{\tau, |u|<|\tau|^{\frac{1}{2}}}\frac{|\tau|^{\frac{1}{4}}}{|\tau|} \left(\int_{|v|< |\tau|^{\frac{1}{2}}}
\frac{|u+v|}{|u-v|}d v\right)^{\frac{1}{2}}.
\end{align*}
After a change of variables this is reduced to $\tau=1$, and $A$ is bounded.
By Sobolev's estimate at an angle, we have 
\begin{align*}
 \big\|{|\nabla_y|^{\frac12}}{|\nabla_x|}^{-\frac 12}f\big\|_{L^2(dtd(x-y)d(x+y))} 
 \lesssim   \| |\nabla_{x}|^{\frac12}|\nabla_{y}|^{\frac12} f\|_{L^{2}(dt)L^{\frac65}(d(x-y))L^{2}(d(x+y))}.
\end{align*}

Similarly, by Cauchy-Schwarz, we have
\bgd\ad&\mathrm{RHS} \eqref{8'}
\lesssim \|\frac{|\tau|^{\frac{1}{4}}}{|\tau|}\int_{|\xi-\eta|< |\tau|^{\frac{1}{2}}} |\mathcal{F}f| d(\xi-\eta)\|_{L^2(d \tau d(\xi+\eta))}\\
&\lesssim A \big\|\dt{|\nabla_{x-y}|^{-\frac12}}f\big\|_{L^2(dt)d(x-y)d(x+y)}
\ead\ed where
\bgd\ad A= \sup_{\tau, \xi+\eta}\frac{1}{|\tau|} \left(\int_{|\xi-\eta|< |\tau|^{\frac{1}{2}}}
|\xi-\eta|\, d(\xi-\eta)\right)^{\frac{1}{2}},
\ead\ed which is bounded.
By Sobolev's inequality, we have 
\begin{align*}
\big\|\dt{|\nabla_{x-y}|^{-\frac12}}f\big\|_{L^2(dt)d(x-y)d(x+y)}
 \lesssim   \| \dt{|\nabla_{x-y}|^{\frac12}} f\|_{L^{2}(dt)L^{\frac65}(d(x-y))L^{2}(d(x+y))}.
\end{align*}
\bigskip

Case 2:  $ |\xi|+ |\eta| > 2 |\tau|^{\frac{1}{2}}$ 

In this case, since we are assuming that $|\tau|^{\frac12}> 10(1+|\xi+\eta|)$.
we have $|\xi-\eta|> |\xi+\eta|$, and also $|\xi-\eta|> |\tau|^{\frac{1}{2}}$. As before, it suffices to bound the right hand side of \eqref{8'}. By Cauchy-Schwarz
\bgd\ad&\mathrm{RHS} \eqref{8'}
\lesssim \|\int_{2|\xi-\eta|> |\xi+\eta|+ |\tau|^{\frac{1}{2}}} |\frac{|\tau|^{\frac{1}{4}}|\mathcal{F}f|}{|\xi-\eta|^2} d(\xi-\eta)\|_{L^2(d \tau d(\xi+\eta))}\\
&\lesssim  A \big\|{|\nabla_y|^{\frac12}}{|\nabla_x|}^{-\frac 12}f\big\|_{L^2(dtd(x-y)d(x+y))} \\
& \lesssim  A \| |\nabla_{x}|^{\frac12}|\nabla_{y}|^{\frac12} f\|_{L^{2}(dt)L^{\frac65}(d(x-y))L^{2}(d(x+y))}
\ead\ed where we used  Sobolev's estimate at an angle in the last inequality. In this case,
\bgd\ad
A^2=\sup_{\xi+\eta, \tau}\int_{2|\xi-\eta|> |\xi+\eta|+ |\tau|^{\frac{1}{2}}} \frac{|\tau|^{\frac12}}{|\xi-\eta|^4}\frac{|\xi|}{|\eta|} d(\xi-\eta).
\ead\ed Again we scale to $|\tau|^{\frac{1}{2}}+ |\xi+\eta|=1$ and have to estimate
\begin{align*}
\int_{|v|>1} \frac{1}{|v|^4}\frac{|u+v|}{|u-v|} dv.
\end{align*}
This is bounded uniformly in $|u|<1$. 

Similarly, since we are assuming $|\xi-\eta|\ge 20N$,  by Cauchy-Schwarz
\begin{align*}
&\mathrm{RHS} \eqref{8'}
\lesssim \|\int_{|\xi-\eta|\ge 20N} \frac{|\tau|^{\frac{1}{4}}|\mathcal{F}f|}{|\xi-\eta|^2} d(\xi-\eta)\|_{L^2(d \tau d(\xi+\eta))}\\
&\lesssim  N^{-\frac12} \big\|\dt f\big\|_{L^2(dtd(x-y)d(x+y))}.
\end{align*}

Case 4 :$ \frac{1}{2}(|\xi|+|\eta|)<|\tau |^{\frac{1}{2}} < 2 (|\xi|+|\eta|)$, 
In this case, since we are assuming  that $|\tau|^{\frac12}> 10(1+|\xi+\eta|)$,
this implies $|\tau |^{\frac{1}{2}} \sim |\xi|\sim |\eta|\sim |\xi-\eta|$.
 
We shall use the decomposition $u=\sum_{k=0}^\infty P_{|\tau|\sim 2^{k}}u$, and  the square function estimate
\bgd\ad\|\dt u\|&_{L^{\infty}(d(x-y)) L^2(d(x+y)dt)} \\ 
&\sim  \|\Big(\sum_{k=0}^\infty |P_{|\tau|\sim 2^{k}} \dt u|^2\Big)^{\frac12}\|_{L^{\infty}(d(x-y)) L^2(d(x+y)dt)} \\
&\lesssim\Big(\sum_{k=0}^\infty 2^{2k} \|P_{|\tau|\sim 2^{k}}u\|_{{L^{\infty}(d(x-y)) L^2(d(x+y)dt)}}^2\Big)^{\frac12}.
\ead\ed
For each fixed dyadic piece $P_{|\tau|\sim 2^{k}}u$, in the current case, we have 
$$P_{|\tau|\sim 2^{k}}u= P_{|\tau|\sim 2^{k}}P_{|\xi|\sim 2^{k/2}}P_{|\eta|\sim 2^{k/2}}u,$$
which implies 
 \begin{align*}
\|\big|P_{|\tau| \sim 2^k}u\|_{collapsing} &\lesssim
\|\dx P_{|\xi|\sim 2^{k/2}}P_{|\eta|\sim 2^{k/2}} u\|_{L^2(dt)L^6(d(x-y)) L^2(d(x+y))}\\ 
& \lesssim
2^{-k/2}\|\dx\dy P_{|\eta|\sim 2^{k/2}} u\|_{L^2(dt)L^6(d(x-y)) L^2(d(x+y))}
\end{align*}
where we used Bernstein's  inequality in rotated coordinates twice, see e.g., Lemma 3.1 in \cite{1}.
Thus,
\bgd\ad&\|\dt u\|_{L^{\infty}(d(x-y)) L^2(d(x+y)dt)} \\
&\lesssim
\Big(\sum_{k=0}^\infty \|\dx\dy P_{|\eta|\sim 2^{k/2}} u\|_{L^2(dt)L^6(d(x-y)) L^2(d(x+y))}^2\Big)^{\frac12}
\\
&\lesssim
\left(\sum_{k=0}^\infty 
\| \dx\dy P_{|\eta| \sim 2^{k/2}}f\|^2_{L^{2}(dt)L^{\frac65}(d(x-y))L^{2}(d(x+y))}\right)^{\frac{1}{2}}\\
& \lesssim \|\dx\dy f\|_{L^{2}(dt)L^{\frac65}(d(x-y))L^{2}(d(x+y))}
\ead\ed where we used Strichartz (Theorem~\ref{2.1'}) in the second line, and square function estimates in $y$ the the last line.

Similarly, for each fixed dyadic piece $P_{|\tau|\sim 2^{k}} \dt u$, in this case since $|\xi-\eta|\sim 2^{\frac k2}$, by Bernstein's  inequality, we have 
\bgd\ad\|\big|&P_{|\tau| \sim 2^k}\dt u\|_{collapsing}\\
 &\lesssim
\|\dt\dyy  P_{|\xi-\eta|\sim 2^{\frac k2}} u\|_{L^2(dt)L^6(d(x-y)) L^2(d(x+y))}\\ 
&\lesssim
\|\dt\dyy  P_{|\xi-\eta|\sim 2^{\frac k2}} f\|_{L^2(dt)L^{\frac65}(d(x-y)) L^2(d(x+y))},
\ead\ed where in the last inequality we used Theorem~\ref{Spsi2}.
Thus,
\bgd\ad &\|\dt u\|_{L^{\infty}(d(x-y)) L^2(d(x+y)dt)} \\
&\lesssim
\left(\sum_{k=0}^\infty 
\| \dt\dyy P_{|\xi-\eta|\sim 2^{\frac k2}}  f\|^2_{L^{2}(dt)L^{\frac65}(d(x-y))L^{2}(d(x+y))}\right)^{\frac{1}{2}}\\
& \lesssim \|\dt\dyy   f\|_{L^{2}(dt)L^{\frac65}(d(x-y))L^{2}(d(x+y))}
\ead\ed where we used square function estimates in $x-y$ the the last line.
\end{proof}

\bibliographystyle{abbrv}
\bibliography{refs}

\end{document}